\RequirePackage{fix-cm}
\documentclass[smallextended]{svjour3}       
\smartqed  
\usepackage{graphicx}
\usepackage{amsmath}
\usepackage{amssymb}
\usepackage{latexsym}
\usepackage{mathrsfs}
\usepackage{colortbl}
\usepackage{amscd}
\usepackage{tikz-cd}
\usepackage{comment}

 \journalname{Japan Journal of Industrial and Applied Mathematics}

\def\div{\mathop{\mathrm{div}}}
\def\diam{\mathop{\mathrm{diam}}}
\def\diag{\mathop{\mathrm{diag}}}
\def\>{\textgreater}
\def\<{\textless}

\def\conv{\mathop{\mathrm{conv}}}
\def\esssup{\mathop{\mathrm{ess.sup}}}

\spnewtheorem{thr}{Theorem}{\bf}{\it}
\spnewtheorem{coro}{Corollary}{\bf}{\it}
\spnewtheorem{defi}{Definition}{\bf}{\it}
\spnewtheorem{lem}{Lemma}{\bf}{\it}
\spnewtheorem{prop}{Proposition}{\bf}{\it}
\spnewtheorem{ass}{Assumption}{\bf}{\it}
\spnewtheorem{cond}{Condition}{\bf}{\it}
\spnewtheorem{Rem}{Remark}{\it}{\it}
\spnewtheorem{Ex}{Example}{\it}{\it}
\spnewtheorem{Note}{Note}{\it}{\it}
\spnewtheorem*{ex*}{Example:}{\it}{\it}
\spnewtheorem*{pf*}{Proof}{\bf}{\rm}
\spnewtheorem*{rem*}{Remark:}{\it}{\it}
\spnewtheorem*{note*}{Note:}{\it}{\it}

\spnewtheorem*{lem1*}{Lemma 1}{\bf}{\rm}
\spnewtheorem*{lem3*}{Lemma 3}{\bf}{\rm}

\newcounter{sone}
\newcounter{stwo}
\newcounter{sthree}
\newcounter{sfour}
\newcounter{sfive}
\newcounter{ssix}
\newcounter{lone}
\newcounter{ltwo}
\newcounter{lthree}
\newcounter{lfour}
\newcounter{lfive}
\newcounter{lsix}
\makeatletter
    
    \@addtoreset{equation}{section}
  \makeatother
  
\begin{document}

\title{General theory of interpolation error estimates on anisotropic meshes
}

\titlerunning{Interpolation error estimates on anisotropic meshes}        

\author{Hiroki Ishizaka \and Kenta Kobayashi \and Takuya Tsuchiya 
}


\institute{Hiroki Ishizaka \at
              Graduate School of Science and Engineering, Ehime University, Matsuyama, Japan \\
              \email{h.ishizaka005@gmail.com}           
           \and
           Kenta Kobayashi\at
              Graduate School of Business Administration, Hitotsubashi University, Kunitachi, Japan \\
            \email{kenta.k@r.hit-u.ac.jp}
            \and
            Takuya Tsuchiya \at 
             Graduate School of Science and Engineering, Ehime University, Matsuyama, Japan \\
              \email{tsuchiya@math.sci.ehime-u.ac.jp}  
}

\date{Received: date / Accepted: date}

\maketitle

\begin{abstract}
We propose a general theory of estimating interpolation error for  smooth functions in two and three dimensions. In our theory, the error of interpolation is bound in terms of the diameter of a simplex and a geometric parameter. In the two-dimensional case, our geometric parameter is equivalent to the circumradius of a triangle. In the three-dimensional case, our geometric parameter also represents the flatness of a tetrahedron. Through the introduction of the geometric  parameter, the error estimates newly obtained can be applied to cases that violate the maximum-angle condition.

\end{abstract}

\section{Introduction}
\label{intro}
It is challenging to construct accurate and efficient finite element schemes for solving partial differential equations in various domains. Estimations of interpolation error are important in terms of ensuring the validity of schemes and their accuracy sometimes depends on geometric conditions of meshes of the domain. Many studies have imposed the condition of shape regularity to a family of meshes \cite{Bra07,BreSco08,Cia02,ErnGue04,Ran17}; i.e., triangles or tetrahedra cannot be too flat in a shape-regular family of triangulations. 

In \cite{BabAzi76}, the shape regularity condition was relaxed to the maximum-angle condition, which refers to the maximum angle of each triangle in meshes being smaller than a constant $\< \pi$. A family of triangulations under the maximum-angle condition allows the use of anisotropic finite element meshes. Anisotropic meshes have different mesh sizes in different directions, and the shape regularity assumption on triangulations is no longer valid on these meshes.

The question arises whether the maximum-angle condition can be relaxed further. The answer was given by \cite{HanKorKri12,KobTsu14,KobTsu15,KobTsu20}; i.e., it is known that the maximum-angle condition is not necessarily needed to obtain error estimates. 

The present paper proposes a general theory of interpolation error estimates for smooth functions that can be applied to, for example, Lagrange, Hermite, and Crouzeix--Raviart interpolations. For a $d$-simplex $T$, we introduce a new geometric parameter $H_T$ in Section \ref{paraH} and the error of interpolations is bounded in terms of the diameter $h_T$ of $T$ and $H_T$. We emphasize that we do not impose the shape regularity condition and the maximum-angle condition for the mesh partition. 

Using the new parameter $H_T$, we also propose error estimates for the Raviart--Thomas interpolation. The Raviart--Thomas interpolation error estimates on anisotropic meshes play an important role in first-order Crouzeix--Raviart finite element analysis. In \cite{AcoDur99}, the interpolation error analysis in the lowest-order case was given under the maximum-angle condition for triangles and tetrahedra. In \cite{AcoApe10}, the authors extended the results to the Raviart--Thomas interpolation with any order in two- and three-dimensional cases. 

Meanwhile, in \cite{KobTsu18}, the lowest-order Raviart--Thomas interpolation error analysis under a condition weaker than the maximum-angle condition was introduced in the two-dimensional case. The analysis was based on the technique of Babu\v{s}ka and Aziz \cite{BabAzi76}. The technique requires a Poincar\'e-like inequality on reference elements.  However, it is not easy to deduce the inequality in the three-dimensional case. To overcome this difficulty, we use the component-wise stability estimates of the Raviart--Thomas interpolation in reference elements introduced in \cite{AcoApe10}. We consequently have the Raviart--Thomas interpolation error estimates of any order in two- and three-dimensional cases under the relaxed mesh condition.

The remainder of the paper is organized as follows. Section 2 introduces notations and basic concepts of the Raviart--Thomas finite element. Section 3 introduces standard positions and the new geometric parameter. Further, we propose affine mappings and Piola transformations on standard positions and present the finite element generation. Section 4 proves interpolation error estimates of smooth functions that can be applied to, for example, Lagrange, Hermite, and Crouzeix--Raviart interpolations. Section 5 proves the Raviart--Thomas interpolation error estimate. Our main theorems are presented as Theorem \ref{int=thr2} and Theorem \ref{rt=thr3}.

\section{Preliminaries} \subsection{Function Spaces}
Let $d=2,3$. Let $\mathbb{N}_0$ denote the set of non-negative integers. Let $\beta := (\beta_1,\ldots,\beta_d)^T \in \mathbb{N}_0^d$ be a multi-index. For the multi-index $\beta$, let
\begin{align*}
\displaystyle
\partial^{\beta} := \left( \frac{\partial}{\partial x_1} \right)^{\beta_1} \ldots  \left( \frac{\partial}{\partial x_d} \right)^{\beta_d} = \frac{\partial^{|\beta|}}{\partial x_1^{\beta_1} \ldots \partial x_d^{\beta_d}} \quad \text{with} \quad |\beta| := \beta_1 + \ldots + \beta_d.
\end{align*}
Let $\Omega$ be an open domain of $\mathbb{R}^d$.  Let $\ell$ be a nonnegative integer and $p \in \mathbb{R}$ with $1 \leq p \leq \infty$. We define the Sobolev space
\begin{align*}
\displaystyle
W^{\ell,p}(\Omega) := \left \{ \varphi \in L^p(\Omega); \ \partial^{\beta} \varphi \in L^p(\Omega), \ 0 \leq  |\beta| \leq \ell \right\},
\end{align*}
equipped with the norms
\begin{align*}
\displaystyle
\| \varphi \|_{W^{\ell,p}(\Omega)} &:= \left (  \sum_{0 \leq |\beta| \leq \ell} \| \partial^{\beta} \varphi \|^p_{L^p(\Omega)} \right)^{1/p} \ \text{if $1 \leq p \< \infty$}, \\
\| \varphi \|_{W^{\ell,\infty}(\Omega)} &:= \max_{0 \leq |\beta| \leq \ell} \left( \esssup_{x \in \Omega} |\partial^{\beta} \varphi (x)| \right).
\end{align*}
We use the semi-norms
\begin{align*}
\displaystyle
| \varphi |_{W^{\ell,p}(\Omega)} &:= \left (  \sum_{|\beta| = \ell} \| \partial^{\beta} \varphi \|^p_{L^p(\Omega)} \right)^{1/p} \ \text{if $1 \leq p \< \infty$}, \\
| \varphi |_{W^{\ell,\infty}(\Omega)} &:= \max_{|\beta| =\ell} \left( \esssup_{x \in \Omega} |\partial^{\beta} \varphi (x)| \right).
\end{align*}
If $p=2$, we use the notation
\begin{align*}
\displaystyle
H^{\ell}(\Omega) := W^{\ell,2}(\Omega).
\end{align*}
We set $L^2(\Omega) := H^{0}(\Omega)$. The space $H^{\ell}(\Omega)$ is a Hilbert space equipped with the scalar product
\begin{align*}
\displaystyle
(\varphi,\psi)_{H^{\ell}(\Omega)} := \sum_{|\beta| \leq \ell} (\partial^{\beta} \varphi , \partial^{\beta} \psi)_{L^2(\Omega)},
\end{align*}
where $(_\cdot, _\cdot)_{L^2(\Omega)}$ denotes the $L^2$-inner product, which leads to the norm and semi-norm
\begin{align*}
\displaystyle
\| \varphi \|_{H^{\ell}(\Omega)} := \left( \sum_{|\beta| \leq \ell} \| \partial^{\beta} \varphi \|^2_{L^2(\Omega)}  \right)^{1/2}, \ | \varphi |_{H^{\ell}(\Omega)} := \left( \sum_{|\beta| = \ell} \| \partial^{\beta} \varphi \|^2_{L^2(\Omega)}  \right)^{1/2}.
\end{align*}

The dual space of $W^{\ell,p}(\Omega)$ is defined $\mathcal{L}(W^{\ell,p}(\Omega);\mathbb{R})$ and denoted by $W^{\ell,p}(\Omega)^{\prime}$. $W^{\ell,p}(\Omega)^{\prime}$ is a Banach space with norm
\begin{align*}
\displaystyle
\| \chi \|_{W^{\ell,p}(\Omega)^{\prime}} := \sup_{v \in W^{\ell,p}(\Omega)} \frac{|\chi(v)|}{\| v \|_{W^{\ell,p}(\Omega)}} \quad \forall \chi \in W^{\ell,p}(\Omega)^{\prime}.
 \end{align*}

For any $v = (v_1,\ldots,v_d)^T \in W^{\ell,p}(\Omega)^d$, the norm is defined by
\begin{align*}
\displaystyle
\| v\|_{W^{\ell,p}(\Omega)^d} := \left( \sum_{i=1}^d \| v_i \|^2_{W^{\ell,p}(\Omega)} \right)^{1/2}.
\end{align*}

We introduce the function space
\begin{align*}
\displaystyle
H(\div;\Omega) := \left\{ v \in L^2(\Omega)^d ; \ \div v \in L^2(\Omega) \right\},
\end{align*}
with the norm
\begin{align*}
\displaystyle
\| v \|_{H(\div;\Omega)} := \left ( \| v \|_{L^2(\Omega)^d}^2 + \| \div v \|^2  \right)^{1/2}.
\end{align*}

Let $A$ be a $d \times d$ matrix, and $\| A \|_2$ denote an operator norm as
\begin{align*}
\displaystyle
\| A \|_2 := \sup_{0 \neq x \in \mathbb{R}^d} \frac{|A x|}{|x|},
\end{align*}
where $| x | := ( \sum_{i=1}^d |x_i|^2)^{1/2}$ for $x \in \mathbb{R}^d$.

\subsection{Raviart--Thomas Finite Element on Simplices}
For any $k \in \mathbb{N}_0$, let $\mathcal{P}^k$ be the space of polynomials with degree at most $k$. $\mathcal{P}^k(D)$ is spanned by the restriction to $D$ of polynomials in $\mathcal{P}^k$, where $D$ is a closed domain. Let $T$ be a $d$-simplex. 
\begin{defi} \label{RT=defi}
The local Raviart--Thomas polynomial space of order $k \in \mathbb{N}_0$ is defined by
\begin{align}
\displaystyle
RT^k(T) := \mathcal{P}^k(T)^d + x \mathcal{P}^k(T), \quad x \in \mathbb{R}^d. \label{pre4}
\end{align}
For $v \in RT^k(T)$, the local degrees of freedom are given as
\begin{align}
\displaystyle
& \int_{F_{i}} v \cdot n_{F_i} p_k ds, \quad \forall p_k \in \mathcal{P}^k(F_i), \quad F_i \subset \partial T, \label{pre5} \\
 &\int_T v \cdot q_{k-1} dx, \quad \forall q_{k-1} \in \mathcal{P}^{k-1}(T)^d. \label{pre6}
\end{align}
Here, $n_{F_i}$ denotes the outer unit normal vector of $T$ on the face $F_i$. Note that for $k=0$, local degrees of freedom of type \eqref{pre6} are violated.
\end{defi}

For the simplicial Raviart--Thomas element in $\mathbb{R}^d$, it holds that
\begin{align}
\displaystyle
\dim RT^k(T) =
\begin{cases}
(k+1)(k+3) \quad \text{if $d=2$}, \\
\frac{1}{2}(k+1)(k+2)(k+4) \quad \text{if $d=3$}.
\end{cases} \label{pre7}
\end{align}
It is known that the Raviart--Thomas finite element with the local degrees of freedom in Definition \ref{RT=defi} is unisolvent; e.g., see \cite[Proposition 2.3.4]{BofBreFor13}. The triple $\{ T , RT^k , \Sigma \}$ is then a finite element.

We set the domain of the local Raviart--Thomas interpolation as $V^{\div}(T) := H^1(T)^d$; e.g., see also \cite[p. 27]{ErnGue04}.

The local Raviart--Thomas interpolation $I_T^{RT}: V^{\div}(T) \to RT^k(T)$ is then defined as follows. For any $v \in V^{\div}(T)$,
\begin{align}
\displaystyle
\int_F I_T^{RT} v \cdot n_F p_k ds &= \int_F v \cdot n_F p_k ds \quad \forall p_k \in \mathcal{P}^k(F), \ F \subset \partial T, \label{pre8}
\end{align}
and if $k \geq 1$,
\begin{align}
\displaystyle
\int_T I_T^{RT} v \cdot q_{k-1} dx &= \int_T v \cdot q_{k-1} dx \quad \forall q_{k-1} \in \mathcal{P}^{k-1}(T)^d. \label{pre9}
\end{align}

Let $\{ \widehat{T} , \widehat{P} , \widehat{\Sigma} \}$ with $\widehat{P} := {RT^k}(\widehat{T})$ be the Raviart--Thomas finite element. Let $ \widehat{\Phi}: \mathbb{R}^d \to \mathbb{R}^d$, ${x} := \widehat{\Phi}(\hat{x}) := \widehat{A} \hat{x} + \hat{b}$  be an affine mapping such that ${T} = \widehat{\Phi} (\widehat{T})$ with a regular matrix  $\widehat{A} \in \mathbb{R}^{d \times d}$ and $\hat{b} \in \mathbb{R}^d$. The Piola transformation $\widehat{\Psi}: L^2(\widehat{T})^d \to L^2({T})^d$ is defined by
\begin{align*}
\displaystyle
\widehat{\Psi}: L^2(\widehat{T})^d \ni \hat{v}(\hat{x}) &\mapsto {v}({x}) := \widehat{\Psi}(\hat{v})({x}) := \frac{1}{| \det (\widehat{A}) |} \widehat{A} \hat{v}(\hat{x}) \in L^2({T})^d.
\end{align*}

The following lemmata introduce the fundamental properties of the Piola transformation. 
\begin{lem} \label{pre=lem1}
For $\hat{\varphi} \in H^1(\widehat{T})$, $\hat{v} \in H(\widehat{\div};\widehat{T})$, we define $\varphi := \hat{\varphi} \circ \widehat{\Phi}^{-1}$ and $v := \widehat{\Psi} (\hat{v})$. Then, 
\begin{align}
\displaystyle
\int_{T} \div v \varphi dx &= \int_{\widehat{T}} \widehat{\div} \hat{v} \hat{\varphi} d \hat{x}, \notag \\
\int_{T} v \cdot \nabla_x \varphi dx &= \int_{\widehat{T}} \hat{v} \cdot \widehat{\nabla}_{\hat{x}} \hat{\varphi} d \hat{x},  \notag \\
\int_{\partial T} v \cdot n_T \varphi ds &= \int_{\partial \widehat{T}} \hat{v} \cdot \hat{n}_{\widehat{T}} \hat{\varphi} d \hat{s}.  \label{pre10}
\end{align}
Here, $n_T$ and $\hat{n}_{\widehat{T}} $ are respectively the unit outward normal vectors of $T$ and $\widehat{T}$.
\end{lem}

\begin{pf*}
See, for example, \cite[Lemma 3.3]{BofBre08}.
\qed
\end{pf*}

By applying \eqref{pre10}, we can prove the invariance of the Raviart--Thomas interpolation under the Piola transform; e.g., see \cite[Lemma 3.4]{BofBre08}.
\begin{lem} \label{pre=lem2}
For $\hat{v} \in H^1(\widehat{T})^d$, we have
\begin{align*}
\displaystyle
I_{\widehat{T}}^{RT} \hat{v} = \widehat{\Psi}^{-1} I_T^{RT} \widehat{\Psi} \hat{v}.
\end{align*}
That is to say, the diagram
\[
\begin{CD}
     H^1(\widehat{T})^d @>{\widehat{\Psi}}>> H^1({T})^d \\
  @V{I_{\widehat{T}}^{RT}}VV    @VV{I_{{T}}^{RT}}V \\
     RT^k(\widehat{T})   @>{\widehat{\Psi}}>>  RT^k(T)
\end{CD}
 \]
commutes. 
\end{lem}

\section{Standard Positions and Reference Elements} \label{reference}
This section introduces the Jacobian matrix proposed in \cite{KobTsu20} for the three-dimensional case and that proposed in \cite{KobTsu14,KobTsu15,LiuKik18} for the two-dimensional case. 

Let us first define a diagonal matrix $ \widehat{A}^{(d)}$ as
\begin{align}
\displaystyle
\widehat{A}^{(d)} :=  \diag (\alpha_1,\ldots,\alpha_d), \quad \alpha_i \in \mathbb{R}.  \label{mesh1}
\end{align}

\subsection{Two-dimensional case} \label{reference2d}
Let $\widehat{T} \subset \mathbb{R}^2$ be the reference triangle with vertices $\hat{x}_1 := (0,0)^T$, $\hat{x}_2 := (1,0)^T$, and $\hat{x}_3 := (0,1)^T$. 

Let $\widetilde{\mathfrak{T}}^{(2)}$ be the family of triangles
\begin{align*}
\displaystyle
\widetilde{T} = \widehat{A}^{(2)} (\widehat{T}),
\end{align*}
with vertices $\tilde{x}_1 := (0,0)^T$, $\tilde{x}_2 := (\alpha_1,0)^T$, and $\tilde{x}_3 := (0,\alpha_2)^T$.

We next define the regular matrices $\widetilde{A} \in \mathbb{R}^{2 \times 2}$ by
\begin{align}
\displaystyle
\widetilde{A} :=
\begin{pmatrix}
1 & s \\
0 & t \\
\end{pmatrix}, \label{mesh2}
\end{align}
with parameters
\begin{align*}
\displaystyle
s^2 + t^2 = 1, \quad t \> 0.
\end{align*}
For $\widetilde{T} \in \widetilde{\mathfrak{T}}^{(2)}$, let $\mathfrak{T}^{(2)}$ be the family of triangles
\begin{align*}
\displaystyle
T &= \widetilde{A} (\widetilde{T}),
\end{align*}
with vertices $x_1 := (0,0)^T, \ x_2 := (\alpha_1,0)^T, \ x_3 :=(\alpha_2 s , \alpha_2 t)^T$. We then have  $\alpha_1 = |x_1 - x_2| \> 0$, $\alpha_2 = |x_1 - x_3| \> 0$. 

\subsection{Three-dimensional cases} \label{reference3d}
Let $\widehat{T}_1$ and $\widehat{T}_2$ be reference tetrahedrons with the following vertices.
\begin{description}
   \item[(\roman{sone})] $\widehat{T}_1$ has the vertices $\hat{x}_1 := (0,0,0)^T$, $\hat{x}_2 := (1,0,0)^T$, $\hat{x}_3 := (0,1,0)^T$.$\hat{x}_4 := (0,0,1)^T$,
 \item[(\roman{stwo})] $\widehat{T}_2$ has the vertices $\hat{x}_1 := (0,0,0)^T$, $\hat{x}_2 := (1,0,0)^T$, $\hat{x}_3 := (1,1,0)^T$.$\hat{x}_4 := (0,0,1)^T$.
\end{description}

Let $\widetilde{\mathfrak{T}}_i^{(3)}$, $i=1,2$, be the family of triangles
\begin{align*}
\displaystyle
\widetilde{T}_i = \widehat{A}^{(3)} (\widehat{T}_i), \quad i=1,2
\end{align*}
with vertices 
\begin{description}
   \item[(\roman{sone})] $\tilde{x}_1 := (0,0,0)^T$, $\tilde{x}_2 := (\alpha_1,0,0)^T$, $\tilde{x}_3 := (0,\alpha_2,0)^T$, and $\tilde{x}_4 := (0,0,\alpha_3)^T$,
 \item[(\roman{stwo})]  $\tilde{x}_1 := (0,0,0)^T$, $\tilde{x}_2 := (\alpha_1,0,0)^T$, $\tilde{x}_3 := (\alpha_1,\alpha_2,0)^T$, and $\tilde{x}_4 := (0,0,\alpha_3)^T$.
\end{description}

We next define the regular matrices $\widetilde{A}_1, \widetilde{A}_2 \in \mathbb{R}^{3 \times 3}$ by
\begin{align}
\displaystyle
\widetilde{A}_1 :=
\begin{pmatrix}
1 & s_1 & s_{21} \\
0 & t_1  & s_{22}\\
0 & 0  & t_2\\
\end{pmatrix}, \
\widetilde{A}_2 :=
\begin{pmatrix}
1 & - s_1 & s_{21} \\
0 & t_1  & s_{22}\\
0 & 0  & t_2\\
\end{pmatrix} \label{mesh3}
\end{align}
with parameters
\begin{align*}
\displaystyle
\begin{cases}
s_1^2 + t_1^2 = 1, \ s_1 \> 0, \ t_1 \> 0, \ \alpha_2 s_1 \leq \alpha_1 / 2, \\
s_{21}^2 + s_{22}^2 + t_2^2 = 1, \ t_2 \> 0, \ \alpha_3 s_{21} \leq \alpha_1 / 2.
\end{cases}
\end{align*}
For $\widetilde{T}_i \in \widetilde{\mathfrak{T}}_i^{(3)}$, $i=1,2$, let $\mathfrak{T}_i^{(3)}$, $i=1,2$ be the family of triangles
\begin{align*}
\displaystyle
T_i &= \widetilde{A}_i (\widetilde{T}_i), \quad i=1,2
\end{align*}
with vertices
\begin{align*}
\displaystyle
&x_1 := (0,0,0)^T, \ x_2 := (\alpha_1,0,0)^T, \ x_4 := (\alpha_3 s_{21},\alpha_3 s_{22},\alpha_3 t_2)^T, \\
&\begin{cases}
x_3 := (\alpha_2 s_1 , \alpha_2 t_1 , 0)^T \quad \text{for the case (\roman{sone})}, \\
x_3 := (\alpha_1 - \alpha_2 s_1, \alpha_2 t_1,0)^T \quad \text{for the case (\roman{stwo})}.
\end{cases}
\end{align*}
We then have $\alpha_1 = |x_1 - x_2| \> 0$, $\alpha_3 = |x_1 - x_4| \> 0$, and
\begin{align*}
\displaystyle
\alpha_2 =
\begin{cases}
|x_1 - x_3| \> 0  \quad \text{for the case (\roman{sone})}, \\
|x_2 - x_3| \> 0  \quad \text{for the case (\roman{stwo})}.
\end{cases}
\end{align*}

\subsection{Standard Positions} \label{APsubsection33}
In what follows, we impose conditions for $T \in \mathfrak{T}^{(2)}$ in the two-dimensional case and $T \in \mathfrak{T}_1^{(3)} \cup \mathfrak{T}_2^{(3)} =: \mathfrak{T}^{(3)}$ in the three-dimensional case. 

\begin{cond}[Case that $d=2$] \label{cond1}
Let $T \in \mathfrak{T}^{(2)}$ with vertices $x_i$ ($i=1,\ldots,3$) introduced in Section \ref{reference2d}. We assume that $\overline{x_2 x_3}$ is the longest edge of $T$; i.e., $ h_T := |x_2 - x_ 3|$. Recall that $\alpha_1 = |x_1 - x_2|$ and $\alpha_2 = |x_1 - x_3|$. We then assume that $\alpha_2 \leq \alpha_1$. Note that $\alpha_1 = \mathcal{O}(h_T)$. 
\end{cond}

\begin{cond}[Case that $d=3$] \label{cond2}
Let $T \in \mathfrak{T}^{(3)}$ with vertices $x_i$ ($i=1,\ldots,4$) introduced in Section \ref{reference3d}. Let $L_i$ ($1 \leq i \leq 6$) be edges of $T$. We denote by $L_{\min}$  the edge of $T$ with minimum length; i.e., $|L_{\min}| = \min_{1 \leq i \leq 6} |L_i|$. Among the four edges that share an end point with $L_{\min}$, we take the longest edge $L^{({\min})}_{\max}$. Let $x_1$ and $x_2$ be end points of the edge $L^{({\min})}_{\max}$. We thus have 
\begin{align*}
\displaystyle
\alpha_1 = |L^{(\min)}_{\max}| = |x_1 - x_2|.
\end{align*}

Consider cutting $\mathbb{R}^3$ with the plane that contains the midpoint of the edge $L^{(\min)}_{\max}$ and is perpendicular to the vector $x_1 - x_2$. We then have two cases: 
\begin{description}
  \item[(Type \roman{sone})] $x_3$ and $x_4$  belong to the same half-space;
  \item[(Type \roman{stwo})] $x_3$ and $x_4$  belong to different half-spaces.
\end{description}
In each case, we respectively set
\begin{description}
  \item[(Type \roman{sone})] $x_1$ and $x_3$ as the end points of $L_{\min}$, that is $\alpha_2 =  |x_1 - x_3| $;
  \item[(Type \roman{stwo})] $x_2$ and $x_3$ as the end points of $L_{\min}$, that is $\alpha_2 =  |x_2 - x_3| $.
\end{description}
Finally, recall that $\alpha_3 = |x_1 - x_4|$. Note that we implicitly assume that $x_1$ and $x_4$ belong to the same half space. Also note that $\alpha_3 \leq 2 \alpha_1$ and $\alpha_1 = \mathcal{O}(h_T)$, where $h_T$ denotes the diameter of $T$.
\end{cond}

\begin{figure}[htbp]
  \includegraphics[width=8cm]{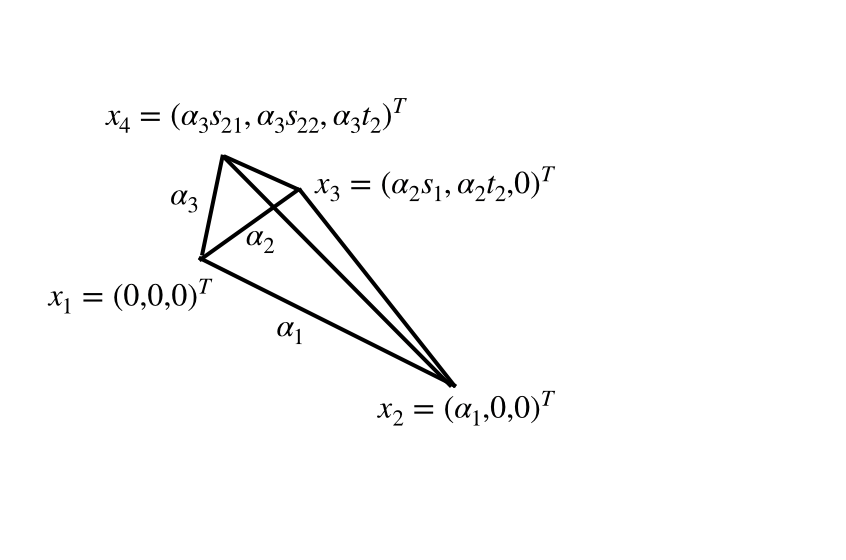}
\caption{Standard position of Type \roman{sone} in $\mathbb{R}^3$}
\label{standardp}
\end{figure}

\begin{figure}[htbp]
 \includegraphics[width=8cm]{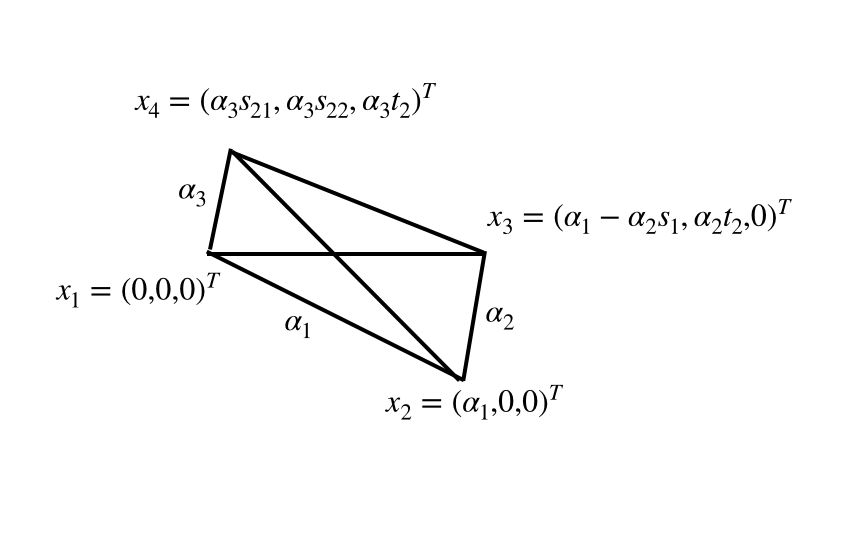}
\caption{Standard position of Type \roman{stwo} in $\mathbb{R}^3$}
\label{standardp2}
\end{figure}

Each $d$-simplex is congruent to the unique $T \in \mathfrak{T}^{(d)}$ satisfying Condition \ref{cond1} or Condition \ref{cond2}. $T$ is therefore called the \textit{standard position} of the $d$-simplex. See Figure \ref{standardp} and \ref{standardp2}.

\subsection{Affine Mappings and Piola Transforms} \label{APsubsection34}
The present paper adopts the following affine mappings and Piola transformations. 

\begin{defi} \label{mesh=defi1}
Let $T \in \mathfrak{T}^{(d)}$ satisfy Condition \ref{cond1} or Condition \ref{cond2}. Let $\widetilde{T}$, and $\widehat{T} \subset \mathbb{R}^d$ be the simplices defined in Sections \ref{reference2d} and \ref{reference3d}. That is to say, 
\begin{align*}
\displaystyle
\widetilde{T} = \widehat{\Phi} (\widehat{T}), \quad {T} = \widetilde{\Phi} (\widetilde{T}) \quad  \text{with} \quad \tilde{x} := \widehat{\Phi}(\hat{x}) := \widehat{A}^{(d)} \hat{x}, \quad x := \widetilde{\Phi}(\tilde{x}) := \widetilde {A} \tilde{x}.
\end{align*}
We define the affine mapping $\Phi: \mathbb{R}^d \to \mathbb{R}^d$ by
\begin{align}
\displaystyle
\Phi := \widetilde{\Phi} \circ \widehat{\Phi}: \mathbb{R}^d \to \mathbb{R}^d, \  {x} :={\Phi}(\hat{x}) := {A} \hat{x}, \quad A := \widetilde {A} \widehat{A}^{(d)}. \label{mesh4}
\end{align}

Let $\widehat{\Psi}: L^2(\widehat{T})^d \to L^2(\widetilde{T})^d$ and $\widetilde{\Psi}: L^2(\widetilde{T})^d \to L^2({T})^d $ be the Piola transformations with respect to $\widehat{A}^{(d)}$ and $\widetilde {A}$, respectively. We define $\Psi  : L^2(\widehat{T})^d  \to L^2(T)^d$ by $\Psi := \widetilde{\Psi} \circ \widehat{\Psi}$, which is the Piola transformation with respect to $A$.
\end{defi}

\subsection{Finite Element Generation on Standard Positions} \label{construction2}
We follow the procedure described in \cite[Section 1.4.1 and 1.2.1]{ErnGue04}.

For the reference element $\widehat{T}$ defined in Sections \ref{reference2d} and \ref{reference3d}, let $\{ \widehat{T} , \widehat{{P}} , \widehat{\Sigma} \}$ be a fixed reference finite element, where $\widehat{{P}} $ is a vector space of functions $\hat{p}: \widehat{T} \to \mathbb{R}^n$ for some positive integer $n$ (typically $n=1$ or $n=d$) and $\widehat{\Sigma}$ is a set of $n_{0}$ linear forms $\{ \hat{\chi}_1 , \ldots , \hat{\chi}_{n_0} \}$ such that
\begin{align*}
\displaystyle
\widehat{{P}} \ni \hat{p} \mapsto (\hat{\chi}_1(\hat{p}) , \ldots , \hat{\chi}_{n_0}(\hat{p}))^T \in \mathbb{R}^{n_0}
\end{align*}
is bijective; i.e., $\widehat{\Sigma}$ is a basis for $\mathcal{L}(\widehat{P};\mathbb{R})$. Further, we denote by $\{ \hat{\theta}_1 , \ldots, \hat{\theta}_{n_0} \}$ in $\widehat{{P}}$ the local ($\mathbb{R}^n$-valued) shape functions such that
\begin{align*}
\displaystyle
\hat{\chi}_i(\hat{\theta}_j) = \delta_{ij}, \quad 1 \leq i,j \leq n_0.
\end{align*}

Let $V(\widehat{T})$ be a normed vector space of functions $\hat{v}: \widehat{T} \to \mathbb{R}^n$ such that $\widehat{P} \subset V(\widehat{T})$ and the linear forms $\{ \hat{\chi}_1 , \ldots , \hat{\chi}_{n_0} \}$ can be extended to $V(\widehat{T})^{\prime}$. The local interpolation operator ${I}_{\widehat{T}}$ is then defined by
\begin{align}
\displaystyle
{I}_{\widehat{T}} : V(\widehat{T}) \ni \hat{v} \mapsto \sum_{i=1}^{n_0} \hat{\chi}_i (\hat{v}) \hat{\theta}_i \in \widehat{{P}}. \label{mesh5}
\end{align}

Let ${\Phi}$, $\widetilde{\Phi}$, and $\widehat{\Phi}$ be the affine mappings defined in \eqref{mesh4}. For $T = \widetilde{\Phi} (\widetilde{T}) = \widetilde{\Phi} \circ \widehat{\Phi} (\widehat{T})$, we first define a Banach space $V(T)$ of $\mathbb{R}^n$-valued functions that is the counterpart of $V(\widehat{T})$ and define a linear bijection mapping by
\begin{align*}
\displaystyle
\psi_T := \psi_{\widehat{T}} \circ \psi_{\widetilde{T}}: V(T) \ni v \mapsto \hat{v} := \psi_T(v) := v \circ \Phi \in   V(\widehat{T}),
\end{align*}
with two linear bijection mappings:
\begin{align*}
\displaystyle
&\psi_{\widetilde{T}} : V({T}) \ni v \mapsto \tilde{v} := \psi_{\widetilde{T}}(v) := v \circ \widetilde{\Phi} \in   V(\widetilde{T}), \\
&\psi_{\widehat{T}} : V(\widetilde{T}) \ni \tilde{v} \mapsto \hat{v} :=  \psi_{\widehat{T}} (\tilde{v}) :=  \tilde{v} \circ \widehat{\Phi} \in V(\widehat{T}).
\end{align*}

Furthermore, the triple $\{ \widetilde{T} , \widetilde{P} , \widetilde{\Sigma}\}$ is defined by
\begin{align*}
\displaystyle
\begin{cases}
\displaystyle
\widetilde{T} = \widehat{\Phi}(\widehat{T}); \\
\displaystyle
 \widetilde{P} = \{ \psi_{\widehat{T}}^{-1}(\hat{p}) ; \ \hat{p} \in \widehat{{P}}\}; \\
\displaystyle
\widetilde{\Sigma} = \{ \{ \tilde{\chi}_{i} \}_{1 \leq i \leq n_0}; \ \tilde{\chi}_{i} = \hat{\chi}_i(\psi_{\widehat{T}}(\tilde{p})), \forall \tilde{p} \in \widetilde{P}, \hat{\chi}_i \in \widehat{\Sigma} \},
\end{cases}
\end{align*}
while the triple $\{ T , {P} , \Sigma\}$ is defined by
\begin{align*}
\displaystyle
\begin{cases}
\displaystyle
T = \widetilde{\Phi}(\widetilde{T}); \\
\displaystyle
{P} = \{ \psi_{\widetilde{T}}^{-1}(\tilde{p}) ; \ \tilde{p} \in \widetilde{{P}}\}; \\
\displaystyle
\Sigma = \{ \{ \chi_{i} \}_{1 \leq i \leq n_0}; \ \chi_{i} = \tilde{\chi}_i(\psi_{\widetilde{T}}(p)), \forall p \in {P}, \tilde{\chi}_i \in \widetilde{\Sigma} \}.
\end{cases}
\end{align*}
$\{ \widetilde{T} , \widetilde{P} , \widetilde{\Sigma}\}$ and $\{ T , {P} , \Sigma\}$ are then finite elements. The local shape functions are $\tilde{\theta}_{i} = \psi_{\widehat{T}}^{-1}(\hat{\theta}_i)$ and $\theta_{i} = \psi_{\widetilde{T}}^{-1}(\tilde{\theta}_i)$, $1 \leq i \leq n_0$, and the associated local interpolation operators are respectively defined by
\begin{align}
\displaystyle
{I}_{\widetilde{T}} : V(\widetilde{T}) \ni \tilde{v} \mapsto {I}_{\widetilde{T}} \tilde{v} &:= \sum_{i=1}^{n_0} \tilde{\chi}_{i}(\tilde{v}) \tilde{\theta}_{i} \in \widetilde{P}, \label{mesh6} \\
{I}_T : V(T) \ni v \mapsto {I}_T v &:= \sum_{i=1}^{n_0} \chi_{i}(v) \theta_{i} \in {P}. \label{mesh7}
\end{align}

\begin{prop} \label{mesh=prop2}
The diagrams
\[
\begin{tikzcd}
V(T) \arrow{r}{\psi_{\widetilde{T}}} \arrow[swap]{d}{{I}_T} & V(\widetilde{T}) \arrow[swap]{d}{{I}_{\widetilde{T}}}  \arrow{r}{\psi_{\widehat{T}}} & V(\widehat{T}) \arrow[swap]{d}{{I}_{\widehat{T}}} \\
{P}  \arrow{r}{\psi_{\widetilde{T}}} & \widetilde{{P}} \arrow{r}{\psi_{\widehat{T}}} & \widehat{{P}}
\end{tikzcd}
\]
commute.
\end{prop}

\begin{pf*}
See \cite[Proposition 1.62]{ErnGue04}.
\qed
\end{pf*}

\begin{Ex}
Let  $\{ \widehat{T} , \widehat{{P}} , \widehat{\Sigma} \}$ be a finite element. 
\begin{enumerate}
 \item For the Lagrange finite element of degree $k$, we set $V(\widehat{T}) := \mathcal{C}^0(\widehat{T})$.
 \item For the Hermite finite element, we set $V(\widehat{T}) := \mathcal{C}^1(\widehat{T})$.
 \item For the Crouzeix--Raviart finite element with $k=1$, we set $V(\widehat{T}) := W^{1,1}(\widehat{T})$. 
\end{enumerate}
\end{Ex}

\subsection{Raviart--Thomas Finite Element on Standard Positions} \label{construction2}
For the reference element $\widehat{T}$ defined in Section \ref{reference2d} and \ref{reference3d}, let $\{ \widehat{T} , RT^k(\widehat{T}) , \widehat{\Sigma} \}$ be the Raviart--Thomas finite element with $k \in \mathbb{N}_0$.  Let ${\Phi}$, $\widetilde{\Phi}$, and $\widehat{\Phi}$ be the affine mappings defined in \eqref{mesh4}. Let ${\Psi}$, $\widetilde{\Psi}$, and $\widehat{\Psi}$ be the Piola transformations defined in Definition \ref{mesh=defi1}.

We then define $\{ \widetilde{T} , RT^k(\widetilde{T}) , \widetilde{\Sigma} \}$ and $\{ {T} , RT^k({T}) , {\Sigma} \}$ by
\begin{align*}
\displaystyle
\begin{cases}
\displaystyle
\widetilde{T} = \widehat{\Phi}(\widehat{T}); \\
\displaystyle
RT^k(\widetilde{T}) = \{ \widehat{\Psi}(\hat{p}) ; \ \hat{p} \in RT^k(\widehat{T}) \}; \\
\displaystyle
\widetilde{\Sigma} = \{ \{ \tilde{\chi}_{i} \}_{1 \leq i \leq n_0}; \ \tilde{\chi}_{i} = \hat{\chi}_i( \widehat{\Psi}^{-1}(\tilde{p})), \forall \tilde{p} \in RT^k(\widetilde{T}), \hat{\chi}_i \in \widehat{\Sigma}  \};
\end{cases}
\end{align*}
and 
\begin{align*}
\displaystyle
\begin{cases}
\displaystyle
T = \widetilde{\Phi}(\widetilde{T}); \\
\displaystyle
RT^k(T) = \{ \widetilde{\Psi}(\tilde{p}) ; \ \tilde{p} \in RT^k(\widetilde{T})\}; \\
\displaystyle
\Sigma = \{ \{ \chi_{i} \}_{1 \leq i \leq n_0}; \ \chi_{i} = \tilde{\chi}_i( \widetilde{\Psi}^{-1}(p)), \forall p \in RT^k(T), \tilde{\chi}_i \in \widetilde{\Sigma} \}.
\end{cases}
\end{align*}
$\{ \widetilde{T} , RT^k(\widetilde{T}) , \widetilde{\Sigma} \}$ and $\{ {T} , RT^k({T}) , {\Sigma} \}$ are then the Raviart--Thomas finite elements. Furthermore, let 
\begin{align}
\displaystyle
I_{\widetilde{T}}^{RT}: V^{\div}(\widetilde{T}) \to RT^k(\widetilde{T}) \label{mesh8}
\end{align}
and
\begin{align}
\displaystyle
I_{{T}}^{RT}: V^{\div}(T) \to RT^k(T) \label{mesh9}
\end{align}
be the associated local Raviart--Thomas interpolation defined in \eqref{pre8} and \eqref{pre9}, respectively.

\subsection{Parameter $H_T$ and Mesh} \label{paraH}
We first propose a new parameter $H_T$.

\begin{defi} \label{mesh=defi3}
 Let $T \in \mathfrak{T}^{(d)}$ satisfy Condition \ref{cond1} or Condition \ref{cond2}. Furthermore, let $\alpha_1, \ldots,\alpha_d$ be defined in Condition \ref{cond1} or Condition \ref{cond2}. We then define the parameter $H_T$ as
\begin{align*}
\displaystyle
H_T := \frac{\prod_{i=1}^d \alpha_i}{|T|} h_T,
\end{align*}
where $h_{T} := \diam(T) = \max_{x_1,x_2 \in T} | x_1 - x_2 |$.
\end{defi}

In the sequel of this paper, the interpolation errors are bounded in terms of $H_T$ and $h_T$. However, the parameters $H_{T_0}$ and $H$ proposed below might be more convenient for the practical computation of finite element methods.

We assume that $\Omega \subset \mathbb{R}^d$ is a bounded polyhedral domain. Let $\mathbb{T}_h = \{ T_0 \}$ be a simplicial mesh of $\overline{\Omega}$, made up of closed $d$-simplices, such as
\begin{align*}
\displaystyle
\overline{\Omega} = \bigcup_{T_0 \in \mathbb{T}_h} T_0,
\end{align*}
with $h := \max_{T_0 \in \mathbb{T}_h} h_{T_0}$, where $ h_{T_0} := \diam(T_0)$. We assume that each face of any $d$-simplex $T_1$ in $\mathbb{T}_h$ is either a subset of the boundary $\partial \Omega$ or a face of another $d$-simplex $T_2$ in  $\mathbb{T}_h$. That is, $\mathbb{T}_h$ is a simplicial mesh of $\overline{\Omega}$ without hanging nodes. 

\begin{defi} \label{mesh=defi4}
 Let any simplex $T_0 \in \mathbb{T}_h$ be transformed into $T$ satisfying Condition \ref{cond1} in the two-dimensional case or Condition \ref{cond2} in the three-dimensional case through appropriate rotation, translation, and mirror imaging. We define the parameter $H_{T_0}$ as
\begin{align*}
\displaystyle
H_{T_0} := \frac{h_{T_0}^2}{|T_0|} \min_{1 \leq i \leq 3} |L_i|  \quad \text{if $d=2$},
\end{align*}
where $L_i$ $(i=1,2,3)$ denotes edges of the triangle $T_0$. Further, we define the parameter $H_{T_0}$ as
\begin{align*}
\displaystyle
H_{T_0} := \frac{h_{T_0}^2}{|T_0|} \min_{1 \leq i , j \leq 6, i \neq j} |L_i| |L_j| \quad \text{if $d=3$},
\end{align*}
where $L_i$ $(i=1,\ldots,6)$ denotes edges of the tetrahedra $T_0$. Here, $|T_0|$ denotes the measure of $T_0$. Furthermore, we set 
\begin{align*}
\displaystyle
H := H(h) := \max_{T_0 \in \mathbb{T}_h} H_{T_0}.
\end{align*}
\end{defi}

We practically impose the following assumption.
 \begin{ass}  \label{mesh=ass1}
 We assume that $\{ \mathbb{T}_h \}_{h \> 0}$ is a sequence of  triangulations of $\Omega$ such that
\begin{align*}
\displaystyle
\lim_{h \to 0} H(h) = 0.
\end{align*}
 \end{ass}
 
 \begin{lem} \label{equilem3}
  Let any simplex $T_0 \in \mathbb{T}_h$ be transformed into the standard position $T$ satisfying Condition \ref{cond1} in the two-dimensional case or Condition \ref{cond2} in the three-dimensional case through appropriate rotation, translation, and mirror imaging. Then, there exist positive constants $c_1$ and $c_2$ such that
\begin{align*}
\displaystyle
c_1 H_{T_0} \leq H_T  \leq c_2 H_{T_0}.
\end{align*}
Furthermore, in the two-dimensional case, $H_{T_0}$ is equivalent to the circumradius of $T_0$. Furthermore, the condition $H_{T_0} / h_{T_0} \< \infty$ implies the semiregularity condition \cite{Kri91}, which is equivalent to the maximum-angle condition.
 \end{lem}
 
  Note that the length of all edges of a simplex and measure of the simplex does not change by the transformation. 

\begin{pf*}
We consider for each dimension, $d=2,3$.
 \begin{description}
   \item[Two-dimensional case]\mbox{}\\
   Let  $L_i$ $(i=1,2,3)$ denote edges of the triangle $T_0$ with $|L_1| \leq |L_2| \leq |L_3|$.  It obviously holds that $\alpha_2 = |L_1|$ and $h_T = |L_3| = h_{T_0}$. Because $\alpha_2 \leq \alpha_1 \< 2 h_T$ and $h_T \< \alpha_1 + \alpha_2 \leq 2 \alpha_1$ for the triangle $\triangle x_1 x_2 x_3$, it holds that
\begin{align*}
\displaystyle
\frac{1}{2} h_{T_0} = \frac{1}{2} h_{T} \< \alpha_1 = |L_2| \< 2 h_T = 2 h_{T_0}. 
\end{align*}
Thus, we have
\begin{align*}
\displaystyle
\frac{1}{2} H_{T_0} = \frac{1}{2} \frac{|L_1|}{|T_0|} h_{T_0}^2 \< H_T = \frac{\alpha_1 \alpha_2}{|T|} h_T \< 2 \frac{|L_1|}{|T_0|} h_{T_0}^2 = 2 H_{T_0}.
\end{align*}
Furthermore, it holds that
\begin{align*}
\displaystyle
2 R_2 = 2 \frac{|L_1| |L_2| |L_3|}{4 |T_0|} \< H_{T_0} = \frac{|L_1|}{|T_0|} h_{T_0}^2 \< 8 \frac{|L_1| |L_2| |L_3|}{4 |T_0|} = 8 R_2,
\end{align*}
where $R_2$ denotes the circumradius of $T_0$.

Set the angle between the segments $\overline{x_1 x_2}$ and $\overline{x_1 x_3}$ by $\theta_{\max}$. Furthermore, we set $t := \sin \theta_{\max}$. By construct of the standard position in the two-dimensional case, the angle $\theta_{\max}$ is the maximum angle of $T$. It is easily proven that the maximum angle condition, i.e., there exists a constant $\delta_{\max} \in (0,\pi)$ such that
\begin{align*}
\displaystyle
\theta_{\max} \leq \delta_{\max}
\end{align*}
is equivalent to the condition
\begin{align*}
\displaystyle
 \frac{H_T}{h_T} = \frac{\alpha_1 \alpha_2}{|T|} = \frac{2}{t} = \frac{2}{\sin \theta_{\max}} \< \infty,
\end{align*}
e.g., see \cite{Kri91}. This implies that the maximum angle condition on $T_0 \in \mathbb{T}_h$ is equivalent to the condition
\begin{align*}
\displaystyle
 \frac{H_{T_0}}{h_{T_0}}  \< \infty.
\end{align*}
   
\item[Three-dimensional case]\mbox{}\\
Let  $L_i$ $(i=1,\ldots,6)$ denote edges of the triangle $T_0$ with $|L_1| \leq |L_2| \leq \cdots \leq |L_6|$.  It obviously holds that $\alpha_2 = |L_1|$ and $h_T = |L_6| = h_{T_0}$. We consider for each type of the standard position introduced in Section \ref{APsubsection33}.

\begin{description}
  \item[(Type \roman{sone})] We set $\alpha_4 := |x_3 - x_4|$, $\alpha_5 := |x_2 - x_4|$, and $\alpha_6 := |x_2- x_3|$.  Because $\alpha_1 = |L^{(\min)}_{\max}| = |x_1 - x_2|$ is the longest edge among the four edges that share an end point with $L_{1}$, it holds that
 \begin{align}
\displaystyle
\alpha_2 \leq \min \{ \alpha_3, \alpha_4, \alpha_6 \} \leq  \max \{ \alpha_3, \alpha_4, \alpha_6 \} \leq \alpha_1. \label{alpha346}
\end{align}
Because $x_1$ and $x_4$ belong to the same half-space for the triangle $\triangle x_1 x_2 x_4$,  it holds that
\begin{align*}
\displaystyle
\begin{cases}
\alpha_3 \leq \alpha_5 \leq \alpha_1 = h_T  \quad \text{or}\\
\alpha_3 \leq \alpha_1 \leq \alpha_5 = h_T.
\end{cases}
\end{align*}
Thus, we have
\begin{align*}
\displaystyle
\begin{cases}
\alpha_3 \leq \alpha_5 \leq \alpha_1 = h_T \quad \text{or}\\
\alpha_3 \leq \alpha_1 \leq h_T  \<  2 \alpha_1, \quad \frac{1}{2} h_T \< \alpha_1 \leq h_T.
\end{cases}
\end{align*}
Because $\alpha_3 \leq \alpha_5$, the length of the edge $L_2$ is equal to the one of $\alpha_3$, $\alpha_4$, or $\alpha_6$. 

Assume that $|L_2| = \alpha_3$. We then have 
\begin{align*}
\displaystyle
\frac{1}{2}  H_{T_0} = \frac{1}{2} \frac{|L_1| |L_2| }{|T_0|} h_{T_0}^2  \< H_{T} = \frac{\alpha_1 \alpha_2 \alpha_3}{|T|} h_T \leq \frac{|L_1| |L_2| }{|T_0|} h_{T_0}^2 = H_{T_0}.
\end{align*}

Assume that $|L_2| = \alpha_4$. We consider the triangle $\triangle x_1 x_3 x_4$. From the assumption, we have $\alpha_2 \leq \alpha_4 \leq \alpha_3$ and $\frac{1}{2} \alpha_3 \< \alpha_4 \leq \alpha_3$. We then obtain
\begin{align*}
\displaystyle
\frac{1}{2}  H_{T_0} = \frac{1}{2} \frac{|L_1| |L_2 |}{|T_0|} h_{T_0}^2  \< H_{T} = \frac{\alpha_1 \alpha_2 \alpha_3}{|T|} h_T \< 2 \frac{|L_1| |L_2| }{|T_0|} h_{T_0}^2 = 2 H_{T_0}.
\end{align*}

Assume that $|L_2| = \alpha_6$. We consider the triangle $\triangle x_1 x_2 x_3$. Because $x_1$ and $x_3$ belong to the same half-space for the triangle $\triangle x_1 x_2 x_3$, it holds that $\alpha_2 \leq \alpha_6 \leq \alpha_1$ and $\frac{1}{2} \alpha_1 \< \alpha_6 \leq \alpha_1$. From \eqref{alpha346}, we have
\begin{align*}
\displaystyle
\frac{1}{2} \alpha_3 \leq \frac{1}{2} \alpha_1 \< \alpha_6 \leq \alpha_1.
\end{align*}
Because $\alpha_6 \leq \alpha_3$, we then obtain
\begin{align*}
\displaystyle
\frac{1}{2}  H_{T_0} = \frac{1}{2} \frac{|L_1| |L_2 |}{|T_0|} h_{T_0}^2  \< H_{T} = \frac{\alpha_1 \alpha_2 \alpha_3}{|T|} h_T \< 2 \frac{|L_1| |L_2|}{|T_0|} h_{T_0}^2 = 2 H_{T_0}.
\end{align*}

\item[(Type \roman{stwo})] We set $\alpha_4 := |x_3 - x_4|$, $\alpha_5 := |x_2 - x_4|$, and $\alpha_6 := |x_1- x_3|$.  Because $\alpha_1 = |L^{(\min)}_{\max}| = |x_1 - x_2|$ is the longest edge among the four edges that share an end point with $L_{1}$, it holds that
 \begin{align}
\displaystyle
\alpha_2 \leq \min \{ \alpha_4, \alpha_5, \alpha_6 \} \leq  \max \{ \alpha_4, \alpha_5, \alpha_6 \} \leq \alpha_1. \label{alpha456}
\end{align}


Because $x_1$ and $x_4$ belong to the same half-space for the triangle $\triangle x_1 x_2 x_4$ and \eqref{alpha456}, it holds that
\begin{align*}
\displaystyle
\alpha_3 \leq \alpha_5 \leq \alpha_1.
\end{align*}
This implies that $\alpha_1 = h_T$. Therefore, the length of the edge $L_2$ is equal to the one of $\alpha_3$, $\alpha_4$, or $\alpha_6$. 

Assume that $|L_2| = \alpha_3$. We then have 
\begin{align*}
\displaystyle
 H_{T} = \frac{\alpha_1 \alpha_2 \alpha_3}{|T|} h_T = \frac{|L_1| |L_2| }{|T_0|} h_{T_0}^2 = H_{T_0}.
\end{align*}

Assume that $|L_2| = \alpha_4$. For the triangle $\triangle x_2 x_3 x_4$, we have
\begin{align*}
\displaystyle
\alpha_2 \leq \alpha_4 \leq \alpha_5 \< 2 \alpha_4.
\end{align*}
Because $\alpha_3 \leq \alpha_5$ and $\alpha_4 \leq \alpha_3$, it holds that
\begin{align*}
\displaystyle
H_{T_0} = \frac{|L_1| |L_2 |}{|T_0|} h_{T_0}^2  \leq H_{T} = \frac{\alpha_1 \alpha_2 \alpha_3}{|T|} h_T \< 2 \frac{|L_1| |L_2| }{|T_0|} h_{T_0}^2 = 2 H_{T_0}.
\end{align*}

Assume that $|L_2| = \alpha_6$. We have $\alpha_1 \< \alpha_2 + \alpha_6 \< 2 \alpha_6$ for the triangle $\triangle x_1 x_2 x_3$.  Therefore, since $\alpha_6 \leq \alpha_3 \leq \alpha_1$, we obtain
\begin{align*}
\displaystyle
H_{T_0} = \frac{|L_1| |L_2 |}{|T_0|} h_{T_0}^2  \leq H_{T} = \frac{\alpha_1 \alpha_2 \alpha_3}{|T|} h_T \< 2 \frac{|L_1| |L_2| }{|T_0|} h_{T_0}^2 = 2 H_{T_0}.
\end{align*}
\end{description}
\end{description}
\qed
\end{pf*}
 
\begin{Rem}
Let $\{ \mathbb{T}_h \}$ be decompositions of a polyhedral into tetrahedra. For any $\mathbb{T}_h \in \{ \mathbb{T}_h \}$, let any simplex $T_0 \in \mathbb{T}_h$ be transformed into the standard position $T$ satisfying Condition \ref{cond2} in the three-dimensional case through appropriate rotation, translation, and mirror imaging. 

We conjecture that the maximum angle condition on $T_0 \in \mathbb{T}_h$ (see \cite{Kri92}) is equivalent to the (semiregular) condition that there exists a positive constant $C_0^{SR}$ such that
\begin{align}
\displaystyle
 \frac{H_{T_0}}{h_{T_0}} \leq C_0^{SR}. \label{remiregular3d}
\end{align}
We here present a geometric condition which satisfies $H_T/h_T \leq C^{SR}$ on the standard positions in the three-dimensional case.

We denote by $\varphi_T$ the angle between the base $\triangle x_1 x_2 x_3$ of $T$ and the segment $\overline{x_1 x_4}$. Let $\theta_1 \leq \theta_2 \leq \theta_3$ be angles of the base of $T$. Assume that there exists a constant $\overline{\theta} \< \pi$ such that 
\begin{align}
\displaystyle
\theta_3 \leq \overline{\theta}, \label{angle1}
\end{align}
and  there exists constants $\overline{\varphi}_1$ and $\overline{\varphi}_2$ with $0\< \overline{\varphi}_1 \leq \overline{\varphi}_2 \< \pi$ such that 
\begin{align}
\displaystyle
\overline{\varphi}_1 \leq \varphi_T \leq\overline{\varphi}_2, \label{angle2}
\end{align}
it then holds that there exists a positive constant $C^{SR}$ such that
\begin{align*}
\displaystyle
\frac{H_T}{h_T} \leq C^{SR}.
\end{align*}

Recall that there are two types' standard position, (Type \roman{sone}) or (Type \roman{stwo}). We denote by $\theta_T$
\begin{description}
  \item[(Type \roman{sone})]the angle between the segments $\overline{x_1 x_2}$ and $\overline{x_1 x_3}$ or
  \item[(Type \roman{stwo})] the angle between the segments $\overline{x_2 x_1}$ and $\overline{x_2 x_3}$.
\end{description}
We set $t_1 := \sin \theta_T$ and $t_2 := \sin \varphi_T$. 

From \eqref{angle1}, it holds that $\frac{\pi}{3} \leq \theta_3$ and
\begin{align*}
\displaystyle
\theta_2,\theta_3 \in \left[ \frac{\pi - \overline{\theta}}{2}, \overline{\theta} \right].
\end{align*}
The  angle $\theta_T$ is not the minimum angle of the base of $T$ because $\alpha_2$ is the minimum edge of $T$. The angle $\theta$ is then either $\theta_2$ or $\theta_3$. Therefore, we have
\begin{align*}
\displaystyle
\theta_T \in \left[ \frac{\pi - \overline{\theta}}{2}, \overline{\theta} \right].
\end{align*}
From this, we have
\begin{align*}
\displaystyle
M_1 := \min \left\{ \sin \frac{\pi - \overline{\theta}}{2} , \sin  \overline{\theta} \right\} \leq \sin \theta_T.
\end{align*}
Furthermore, from \eqref{angle2}, we have
\begin{align*}
\displaystyle
M_2 := \min \left\{ \sin \overline{\varphi}_1 , \sin  \overline{\varphi}_2\right\} \leq \sin \varphi_T.
\end{align*}
We then conclude
\begin{align*}
\displaystyle
\frac{H_T}{h_T} = \frac{\alpha_1 \alpha_2 \alpha_3}{|T|}  = \frac{6}{t_1 t_2} = \frac{6}{ \sin \theta_T \sin \varphi_T} \leq \frac{6}{M_1 M_2} \< \infty,
\end{align*}
where $|T| = \frac{1}{6} \alpha_1 \alpha_2 \alpha_3 t_1 t_2$.

Especially, the standard position $T$ with $\theta_T = \varphi_T = \frac{\pi}{2}$ satisfies both the maximum angle condition and $H_T/h_T \leq C^{SR}$. On other cases, we need further investigation and we leave them for future work.
\end{Rem}

\begin{Rem}
In \cite{KobTsu20}, the projected circumradius $R_{T_0}$ of a tetrahedron $T_0$ is proposed as a geometric parameter for the three-dimensional case. The parameter $H_{T_0}$ that we here propose is much simpler than $R_{T_0}$. We conjecture that $H_{T_0}$ is equivalent to $R_{T_0}$.
\end{Rem}


\section{Interpolation Error Estimates of Smooth Functions} \label{smooth}
This section proposes interpolation error estimates of smooth functions. 

We first give an estimate related to the diagonal matrix \eqref{mesh1} adopting the Babu\v{s}ka--Aziz technique \cite{BabAzi76}.

\begin{lem} \label{int=lem3}
Let $1 \leq p \leq \infty$ and $k \geq 0$. Let $\ell$ be such that $0 \leq \ell \leq k$. Let $\hat{\varphi} \in W^{m,p}(\widehat{T})$ and $\hat{\psi} \in W^{\ell+1,p}(\widehat{T})$. It then holds that, for all $m \in \{ 0,\ldots,\ell+1\}$,
\begin{align}
\displaystyle
\frac{|\tilde{\varphi}|_{W^{m,p}(\widetilde{T})}}{|\tilde{\psi}|_{W^{\ell+1,p}(\widetilde{T})}}
&\leq \max_{1 \leq i \leq d} \{ \alpha_i^{\ell+1-m} \} \left(  \sum_{|\beta| = m} \|  \partial^{\beta} ( \alpha_{1}^{- \beta_1} \cdots \alpha_{d}^{- \beta_d} \hat{\varphi}) \|^p_{L^p(\widehat{T})} \right)^{1/p} \notag\\
&\quad \times  \left( \sum_{|\delta| = \ell+1 - m}  \sum_{|\beta| = m} \| \partial^{\delta} \partial^{\beta} ( \alpha_{1}^{- \beta_1} \cdots \alpha_{d}^{- \beta_d} \hat{\psi} ) \|^p_{L^p(\widehat{T})} \right)^{-1/p},  \label{int1}
\end{align}
with $\tilde{\varphi} := \hat{\varphi} \circ \widehat{\Phi}^{-1}$ and $\tilde{\psi} := \hat{\psi} \circ \widehat{\Phi}^{-1}$. 
\end{lem}

\begin{pf*}
Let  $\beta$, $\gamma$ and $\delta$ be multi-indices with $|\beta| = m$, $|\gamma| = \ell + 1$ and $|\delta| = \ell + 1 - m$.
 
We first have, from $\hat{x}_j = \alpha_j^{-1} \tilde{x}_j$, that
\begin{align*}
\displaystyle
\partial^{\beta} \tilde{\varphi} &=  \alpha_{1}^{- \beta_1} \cdots \alpha_{d}^{- \beta_d} \partial^{\beta} \hat{\varphi}.
\end{align*}
If $1 \leq p \< \infty$, through a change in variable, we obtain
\begin{align*}
\displaystyle
|\tilde{\varphi}|_{W^{m,p}(\widetilde{T})}^p
&= \sum_{|\beta| = m}  \| \partial^{\beta} \tilde{\varphi} \|^p_{L^p(\widetilde{T})} 
= | \det (\widehat{A}^{(d)}) | \sum_{|\beta| = m} \|  \partial^{\beta} ( \alpha_{1}^{- \beta_1} \cdots \alpha_{d}^{- \beta_d} \hat{\varphi}) \|^p_{L^p(\widehat{T})}.
\end{align*}
We similarly have
\begin{align*}
\displaystyle
&|\tilde{\psi}|_{W^{\ell+1,p}(\widetilde{T})}^p \\
&= \sum_{|\gamma| = \ell+1}  \| \partial^{\gamma} \tilde{\psi} \|^p_{L^p(\widetilde{T})} \\
&=  \sum_{|\delta| = \ell+1 - m}  \sum_{|\beta| = m} \| \partial^{\delta} \partial^{\beta} \tilde{\psi} \|^p_{L^p(\widetilde{T})}\\
&=  | \det (\widehat{A}^{(d)}) | \sum_{|\delta| = \ell+1 - m}  \sum_{|\beta| = m} ( \alpha_{1}^{- \delta_1 - \beta_1} \cdots \alpha_{d}^{- \delta_d - \beta_d} )^p \| \partial^{\delta} \partial^{\beta} \hat{\psi} \|^p_{L^p(\widehat{T})}\\
&\geq  | \det (\widehat{A}^{(d)}) | \min_{1 \leq i \leq d} \{ \alpha_i^{- |\delta| p} \}  \sum_{|\delta| = \ell+1 - m}  \sum_{|\beta| = m} \| \partial^{\delta} \partial^{\beta} ( \alpha_{1}^{- \beta_1} \cdots \alpha_{d}^{- \beta_d} \hat{\psi} ) \|^p_{L^p(\widehat{T})}.
\end{align*}

When $p = \infty$, a proof can be made by analogous argument.
\qed
\end{pf*}

We next give estimates relating to the matrix \eqref{mesh2} and \eqref{mesh3}. To this end, we use the fact that if $A^T A$ is a positive definite matrix in $\mathbb{R}^{d \times d}$, the spectral norm of the matrix $A^T A$ is the largest eigenvalue of $A^T A$; i.e.,
\begin{align*}
\displaystyle
\| A \|_2 = \left( \lambda_{\max}(A^T A) \right)^{1/2} = \sigma_{\max}(A),
\end{align*}
where $\lambda_{\max}(A)$ and $\sigma_{\max}(A)$ are respectively the largest eigenvalues and singular values of $A$.

\begin{lem} \label{int=lem4}
Let $1 \leq p \leq \infty$ and $k \geq 0$. Let $\ell$ be such that $0 \leq \ell \leq k$. Let $\hat{\varphi} \in W^{m,p}(\widetilde{T})$ and $\hat{\psi} \in W^{\ell+1,p}(\widetilde{T})$. It then holds that, for all $m \in \{ 0,\ldots,\ell+1\}$,
\begin{align}
\displaystyle
\frac{| {\varphi} |_{W^{m,p}({T})}}{| {\psi} |_{W^{\ell+1,p}({T})}}
&\leq C^{A,d} \left( \frac{H_T}{h_T} \right)^m \frac{| \tilde{\varphi} |_{W^{m,p}(\widetilde{T})}}{| \tilde{\psi} |_{W^{\ell+1,p}(\widetilde{T})}}, \label{int2}
\end{align}
with ${\varphi} := \tilde{\varphi} \circ \widetilde{\Phi}^{-1}$ and ${\psi} := \tilde{\psi} \circ \widetilde{\Phi}^{-1}$. Here, $C^{A,2} := \sqrt{2}^{\ell+1-m} C^{sc}$, and $C^{A,3} := \frac{2^{\ell+1}}{3^{m}} C^{sc}$, where $C^{sc}$ is a constant independent of $T$ and $\widetilde{T}$.
\end{lem}

\begin{pf*}
Using the standard estimates in \cite[Lemma 1.101]{ErnGue04}, we easily get
\begin{align}
\displaystyle
\frac{| {\varphi} |_{W^{m,p}({T})}}{| {\psi} |_{W^{\ell+1,p}({T})}} &\leq C^{sc} \left( \| \widetilde{A} \|_2 \| \widetilde{A}^{-1} \|_2 \right)^m \| \widetilde{A} \|_2^{\ell+1-m} \frac{| \tilde{\varphi} |_{W^{m,p}(\widetilde{T})}}{| \tilde{\psi} |_{W^{\ell+1,p}(\widetilde{T})}}. \label{int3}
\end{align}

\subsubsection*{Two-dimensional case}
Let $\widetilde{A}$ be introduced in \eqref{mesh2}.
From
\begin{align*}
\displaystyle
\widetilde{A}^T \widetilde{A} =
\begin{pmatrix}
1 & s \\
s & 1 \\
\end{pmatrix}, \quad
\widetilde{A}^{-1} \widetilde{A}^{-T} = \frac{1}{t^2}
\begin{pmatrix}
1 & -s \\
-s & 1 \\
\end{pmatrix},
\end{align*}
we have
\begin{align}
\displaystyle
 \| \widetilde{A} \|_2 
 &= \lambda_{\max}(\widetilde{A}^T \widetilde{A})^{1/2} \leq (1+|s|)^{1/2} \leq \sqrt{2}, \label{int4}
\end{align}
and
\begin{align}
\displaystyle
 \| \widetilde{A} \|_2 \| \widetilde{A}^{-1} \|_2
 &=  \lambda_{\max}(\widetilde{A}^T \widetilde{A})^{1/2}  \lambda_{\max} (\widetilde{A}^{-1} \widetilde{A}^{- T} )^{1/2}
 \leq \frac{2}{t}  = \frac{\alpha_1 \alpha_2}{|T|}, \label{int5}
\end{align}
where we used the fact that $|T| = \frac{1}{2} \alpha_1 \alpha_2 t$.

\subsubsection*{Three-dimensional case}
The matrices $\widetilde{A}_1$ and $\widetilde{A}_2$ introduced in \eqref{mesh3} can be decomposed as $\widetilde{A}_1 = \widetilde{M}_0 \widetilde{M}_1$ and $\widetilde{A}_2 = \widetilde{M}_0 \widetilde{M}_2$ with
\begin{align*}
\displaystyle
\widetilde{M}_0 :=
\begin{pmatrix}
1 & 0 & s_{21} \\
0 & 1  & s_{22}\\
0 & 0  & t_2\\
\end{pmatrix}, \
\widetilde{M}_1 :=
\begin{pmatrix}
1 &  s_1 & 0 \\
0 & t_1  & 0\\
0 & 0  & 1\\
\end{pmatrix}, \
\widetilde{M}_2 :=
\begin{pmatrix}
1 & -s_1 & 0 \\
0 & t_1  & 0\\
0 & 0  & 1\\
\end{pmatrix}.
\end{align*}
\noindent
The eigenvalues of $\widetilde{M}_2^T \widetilde{M}_2$ coincide with those of $\widetilde{M}_1^T \widetilde{M}_1$, and we may therefore suppose without loss of generality that we have Case  (\roman{sone}).

We have the inequalities
\begin{align}
\displaystyle
 \| \widetilde{A}_1 \|_2
&= \lambda_{\max}(\widetilde{A}_1^T \widetilde{A}_1)^{1/2}
\leq \lambda_{\max}(\widetilde{M}_0^T \widetilde{M}_0)^{1/2} \lambda_{\max}(\widetilde{M}_1^T \widetilde{M}_1)^{1/2} \notag\\
&\leq \left(1 +  \sqrt{s_{21}^2 + s_{22}^2} \right)^{1/2} (1 + |s_1|)^{1/2} \leq {2}, \label{int6}
\end{align}
and
\begin{align}
\displaystyle
 \| \widetilde{A}_1 \|_2 \| \widetilde{A}_1^{-1} \|_2
&= \lambda_{\max}(\widetilde{A}_1^T \widetilde{A}_1)^{1/2} \lambda_{\max} (\widetilde{A}_1^{-1} \widetilde{A}_1^{- T} )^{1/2} \notag\\
&\leq \frac{\left(1 +  \sqrt{s_{21}^2 + s_{22}^2} \right) (1 + |s_1| )}{t_1 t_2} 
\leq \frac{4}{t_1 t_2} = \frac{2}{3} \frac{\alpha_1 \alpha_2 \alpha_3}{|T|}, \label{int7}
\end{align}
where we used the fact that $|T| = \frac{1}{6} \alpha_1 \alpha_2 \alpha_3 t_1 t_2$.

Therefore, \eqref{int2} follows from \eqref{int3}, \eqref{int4}, and \eqref{int5} if $d=2$ and from \eqref{int3}, \eqref{int6}, and \eqref{int7} if $d=3$.
\qed
\end{pf*}

To give local interpolation error estimates, we use the inequality given in  \cite[Theorem 1.1]{DekLev04} which is a variant of the Bramble--Hilbert lemma; see also \cite{BreSco08,Ver99}.
\begin{thr} \label{int=thr1}
Let ${D} \subset \mathbb{R}^d$ be a bounded convex domain. Let $\varphi \in W^{m,p}({D})$ with $m \in \mathbb{N}$ and $1 \leq p \leq \infty$. There exists a polynomial $\eta \in \mathcal{P}^{m-1}$ such that
\begin{align}
\displaystyle
 | \varphi - \eta |_{ W^{k,p}({D})} \leq C^{BH}(d,m) \diam({D})^{m-k} | \varphi |_{ W^{m,p}({D})}, \quad k=0,1,\ldots,m.  \label{int8}
\end{align}
\end{thr}

\begin{Rem}
In \cite[Lemma 4.3.8]{BreSco08}, the Bramble--Hilbert lemma is given as follows. Let $B$ be a ball in $D \subset \mathbb{R}^d$ such that $D$ is star-shaped with respect to $B$ and its radius $r \> \frac{1}{2} r_{\max}$, where $r_{\max} := \sup \{ r: D$ is star-shaped with respect to a ball of radius $r \}$. Let $\varphi \in W^{m,p}(D)$ with $m \in \mathbb{N}$ and $1 \leq p \leq \infty$. There exists a polynomial $\eta \in \mathcal{P}^{m-1}$ such that
\begin{align*}
\displaystyle
 | \varphi - \eta |_{ W^{k,p}({D})} \leq C^{BH}(d,m,\gamma) \diam({D})^{m-k} | \varphi |_{ W^{m,p}({D})}, \quad k=0,1,\ldots,m.
\end{align*}
Here, $\gamma$ is called the chunkiness parameter of $D$, which is defined by
\begin{align*}
\displaystyle
\gamma := \frac{\diam({D})}{r_{\max}}.
\end{align*}
The main drawback is that the constant $C^{BH}(d,m,\gamma)$ depends on the chunkiness parameter. Meanwhile, the constant $ C^{BH}(d,m)$ of the estimate \eqref{int8} does not depend on the geometric parameter $\gamma$.
\end{Rem}

\begin{Rem}
For general Sobolev spaces $W^{m,p}(\Omega)$, the upper bounds on the constant $C^{BH}(d,m)$ are not given, as far as we know. However, when $p=2$, the following result has been obtained by Verf$\it{\ddot{u}}$rth \cite{Ver99}.

Let ${D} \subset \mathbb{R}^d$ be a bounded convex domain.  Let $\varphi \in H^{m}({D})$ with $m \in \mathbb{N}$.  There exists a polynomial $\eta \in \mathcal{P}^{m-1}$ such that
\begin{align*}
\displaystyle
 | \varphi - \eta |_{ H^{k}({D})} \leq C^{BH}(d,k,m) \diam({D})^{m-k} | \varphi |_{ H^{m}({D})}, \quad k=0,1,\ldots,m-1.
\end{align*}
Verf$\it{\ddot{u}}$rth has given upper bounds on the constants in the estimates such that
\begin{align*}
\displaystyle
 C^{BH}(d,k,m) \leq \pi^{k-m}
\begin{pmatrix}
 d + k -1 \\
 k
\end{pmatrix}
^{1/2} 
\frac{\{ (m-k) ! \}^{1/2}}{\{ \left[ \frac{m-k}{d} \right] !  \}^{d/2}},
\end{align*}
where $[x]$ denotes the largest integer less than or equal to $x$. 

As an example, let us consider the case $d=3$, $k=1$, and $m=2$. We then have
\begin{align*}
\displaystyle
 C^{BH}(3,1,2) \leq \frac{\sqrt{3}}{\pi},
\end{align*}
thus on the standard reference element $\widehat{T}$ introduced in Section \ref{reference3d}, we obtain
\begin{align*}
\displaystyle
 | \hat{\varphi} - \hat{\eta} |_{ H^{1}({\widehat{T}})} \leq  \frac{\sqrt{6}}{\pi} | \hat{\varphi} |_{ H^{2}({\widehat{T}})} \quad \forall \hat{\varphi} \in H^{2}({\widehat{T}}),
\end{align*}
becase $\diam (\widehat{T}) = \sqrt{2}$.

\end{Rem}

From Theorem \ref{int=thr1}, we have the following estimates.
\begin{thr} \label{int=thr2}
Let $\{ \widehat{T} , \widehat{{P}} , \widehat{\Sigma} \}$ be a finite element with normed vector space $V(\widehat{T})$. Let $1 \leq p \leq \infty$ and assume that there exists a nonnegative integer $k$ such that
\begin{align}
\displaystyle
\mathcal{P}^{k} \subset \widehat{{P}} \subset W^{k+1,p}(\widehat{T}) \subset V(\widehat{T}). \label{int9}
\end{align}
Let $\ell$ ($0 \leq \ell \leq k$) be such that $W^{\ell+1,p}(\widehat{T}) \subset V(\widehat{T})$ with continuous embedding. Let $\Phi$ be an affine mapping defined in \eqref{mesh4} and let $I_T$ be the local interpolation operator on $T$ defined in \eqref{mesh7}. It then holds that for arbitrary $m \in \{ 0 , \ldots, \ell+1 \}$,
\begin{align}
\displaystyle
 |\varphi - {I}_T \varphi|_{W^{m,p}(T)}
&\leq C^I \left( \frac{H_{T}}{h_{T}} \right)^m h_{T}^{\ell+1-m} | \varphi |_{W^{\ell+1,p}(T)}, \label{int10}
\end{align}
for any $\varphi \in W^{\ell+1,p}(T)$. Here, $C^I := C^{A,d} C^{BH} \diam(\widehat{T})^{\ell+1-m} \{ (\ell+1)  C^{\mathcal{S}} \}$ is a positive constant independent of $T$, where $C^{\mathcal{S}}$ is the constant appearing in the proof and $C^{BH} := C^{BH}(d,\ell)$ is the constant appearing in Theorem \ref{int=thr1}.
\end{thr}

\begin{pf*}
Let $\hat{\varphi} \in W^{\ell+1,p}(\widehat{T})$. Let ${I}_{\widetilde{T}}$ and ${I}_{\widehat{T}}$ be the local interpolation operator on $\widetilde{T}$ and $\widehat{T}$ defined in \eqref{mesh5} and \eqref{mesh6}, respectively. From \eqref{int2}, we have
\begin{align}
\displaystyle
\frac{| \varphi - {I}_T \varphi |_{W^{m,p}({T})}}{| {\varphi} |_{W^{\ell+1,p}({T})}}
&\leq C^{A,d} \left( \frac{H_T}{h_T} \right)^m \frac{| \tilde{\varphi} - I_{\widetilde{T}} \tilde{\varphi} |_{W^{m,p}(\widetilde{T})}}{| \tilde{\varphi} |_{W^{\ell+1,p}(\widetilde{T})}}. \label{int11}
\end{align}
Because $0 \leq \ell \leq k$, $\mathcal{P}^{\ell} \subset \mathcal{P}^k \subset \widetilde{{P}}$. Therefore, for any $\tilde{\eta} \in \mathcal{P}^{\ell}$, we have $ I_{\widetilde{T}} \tilde{\eta} = \tilde{\eta}$.
This means that $\mathcal{P}^{\ell}$ is invariant under $I_{\widetilde{T}}$. Using the triangle inequality, we have
\begin{align*}
\displaystyle
\frac{| \tilde{\varphi} - I_{\widetilde{T}} \tilde{\varphi} |_{W^{m,p}(\widetilde{T})}}{| \tilde{\varphi} |_{W^{\ell+1,p}(\widetilde{T})}}
\leq \frac{| \tilde{\varphi} -\tilde{\eta} |_{W^{m,p}(\widetilde{T})}}{| \tilde{\varphi} |_{W^{\ell+1,p}(\widetilde{T})}} + \frac{| I_{\widetilde{T}} ( \tilde{\eta} - \tilde{\varphi}) |_{W^{m,p}(\widetilde{T})}}{| \tilde{\varphi} |_{W^{\ell+1,p}(\widetilde{T})}}.
\end{align*}

Let  $\beta$, $\gamma$ and $\delta$ be multi-indices with $|\beta| = m$, $|\gamma| = \ell + 1$ and $|\delta| = \ell + 1 - m$. Then, using the inequality \eqref{int1}, we have
\begin{align*}
\displaystyle
\frac{| \tilde{\varphi} -\tilde{\eta} |_{W^{m,p}(\widetilde{T})}}{| \tilde{\varphi} |_{W^{\ell+1,p}(\widetilde{T})}}
&\leq h_T^{\ell+1-m} \left(  \sum_{|\beta| = m} \|  \partial^{\beta} ( \alpha_{1}^{- \beta_1} \cdots \alpha_{d}^{- \beta_d} (\hat{\varphi} -  \hat{\eta}) ) \|^p_{L^p(\widehat{T})} \right)^{1/p} \notag\\
&\quad \times  \left( \sum_{|\gamma| = \ell+1} \| \partial^{\gamma} ( \alpha_{1}^{- \beta_1} \cdots \alpha_{d}^{- \beta_d} \hat{\varphi} ) \|^p_{L^p(\widehat{T})} \right)^{-1/p}.
\end{align*}
We thus apply Theorem \ref{int=thr1} to obtain
\begin{align*}
\displaystyle
&h_T^{\ell+1-m} \inf_{\hat{\eta} \in \mathcal{P}^{\ell}(\widehat{T})} \left(  \sum_{|\beta| = m} \|  \partial^{\beta} ( \alpha_{1}^{- \beta_1} \cdots \alpha_{d}^{- \beta_d} (\hat{\varphi} -  \hat{\eta}) ) \|^p_{L^p(\widehat{T})} \right)^{1/p} \\
&\quad \times  \left( \sum_{|\gamma| = \ell+1} \| \partial^{\gamma} ( \alpha_{1}^{- \beta_1} \cdots \alpha_{d}^{- \beta_d} \hat{\varphi} ) \|^p_{L^p(\widehat{T})} \right)^{-1/p} \\
&\leq h_T^{\ell+1-m} C^{BH}(d,\ell) \diam(\widehat{T})^{\ell+1-m}.
\end{align*}

For any $\hat{v} \in V(\widehat{T})$, it holds that
\begin{align*}
\displaystyle
\alpha_{1}^{- \beta_1} \cdots \alpha_{d}^{- \beta_d} I_{\widehat{T}} \hat{v}
= I_{\widehat{T}} (\alpha_{1}^{- \beta_1} \cdots \alpha_{d}^{- \beta_d} \hat{v}),
\end{align*}
and
\begin{align*}
\displaystyle
 &\left(  \sum_{|\beta| = m} \|  \partial^{\beta} ( \alpha_{1}^{- \beta_1} \cdots \alpha_{d}^{- \beta_d} I_{\widehat{T}} \hat{v} )\|^p_{L^p(\widehat{T})} \right)^{1/p} \\
 &\quad = | I_{\widehat{T}} (\alpha_{1}^{- \beta_1} \cdots \alpha_{d}^{- \beta_d} \hat{v}) |_{W^{m,p}(\widehat{T})} \\
 &\quad \leq \sum_{i=1}^{n_0} | \hat{\chi}_i (\alpha_{1}^{- \beta_1} \cdots \alpha_{d}^{- \beta_d} \hat{v}) | |\hat{\theta}_i|_{W^{m,p}(\widehat{T})} \\
 &\quad \leq \sum_{i=1}^{n_0}\|  \hat{\chi}_i \|_{W^{\ell,p}(\Omega)^{\prime}} \| \alpha_{1}^{- \beta_1} \cdots \alpha_{d}^{- \beta_d} \hat{v} \|_{W^{\ell+1,p}(\widehat{T})}  |\hat{\theta}_i|_{W^{m,p}(\widehat{T})} \\
 &\quad \leq C^{S} \| \alpha_{1}^{- \beta_1} \cdots \alpha_{d}^{- \beta_d} \hat{v} \|_{W^{\ell+1,p}(\widehat{T})},
\end{align*}
where $C^S := n_0 \max_{1 \leq i \leq n_0} \|  \hat{\chi}_i \|_{W^{\ell,p}(\Omega)^{\prime}}   |\hat{\theta}_i|_{W^{m,p}(\widehat{T})}$.

Therefore, from the above inequality and inequality \eqref{int1}, we have
\begin{align*}
\displaystyle
\frac{| I_{\widetilde{T}} (\tilde{\varphi} -\tilde{\eta}) |_{W^{m,p}(\widetilde{T})}}{| \tilde{\varphi} |_{W^{\ell+1,p}(\widetilde{T})}}
&\leq h_T^{\ell+1-m} \left(  \sum_{|\beta| = m} \|  \partial^{\beta} ( \alpha_{1}^{- \beta_1} \cdots \alpha_{d}^{- \beta_d}I_{\widehat{T}} (\hat{\varphi} -  \hat{\eta}) ) \|^p_{L^p(\widehat{T})} \right)^{1/p} \notag\\
&\quad \times  \left( \sum_{|\gamma| = \ell+1} \| \partial^{\gamma} ( \alpha_{1}^{- \beta_1} \cdots \alpha_{d}^{- \beta_d} \hat{\varphi} ) \|^p_{L^p(\widehat{T})} \right)^{-1/p} \\
&\leq C^S h_T^{\ell+1-m} \left( \sum_{j=1}^{\ell+1} | \alpha_{1}^{- \beta_1} \cdots \alpha_{d}^{- \beta_d} (\hat{\varphi} -  \hat{\eta}) |^p_{W^{j,p}(\widehat{T})} \right)^{1/p}  \notag\\
&\quad \times  \left( \sum_{|\gamma| = \ell+1} \| \partial^{\gamma} ( \alpha_{1}^{- \beta_1} \cdots \alpha_{d}^{- \beta_d} \hat{\varphi} ) \|^p_{L^p(\widehat{T})} \right)^{-1/p}.
\end{align*}
We apply Theorem \ref{int=thr1} to obtain
\begin{align*}
\displaystyle
& C^S h_T^{\ell+1-m} \inf_{\hat{\eta} \in \mathcal{P}^{\ell}(\widehat{T})}  \left( \sum_{j=1}^{\ell+1} | \alpha_{1}^{- \beta_1} \cdots \alpha_{d}^{- \beta_d} (\hat{\varphi} -  \hat{\eta}) |^p_{W^{j,p}(\widehat{T})} \right)^{1/p} \\
&\quad \times  \left( \sum_{|\gamma| = \ell+1} \| \partial^{\gamma} ( \alpha_{1}^{- \beta_1} \cdots \alpha_{d}^{- \beta_d} \hat{\varphi} ) \|^p_{L^p(\widehat{T})} \right)^{-1/p} \\
&\leq (\ell+1) C^S C^{BH}(d,\ell) \diam(\widehat{T})^{\ell+1-m} h_T^{\ell+1-m}.
\end{align*}
We thus have
\begin{align}
\displaystyle
\frac{| \tilde{\varphi} - I_{\widetilde{T}} \tilde{\varphi} |_{W^{m,p}(\widetilde{T})}}{| \tilde{\varphi} |_{W^{\ell+1,p}(\widetilde{T})}}
&\leq C^{BH} \diam(\widehat{T})^{\ell+1-m} \{ 1+(\ell+1)  C^{\mathcal{S}} \} h_T^{\ell+1-m}.  \label{int12}
\end{align}

We conclude from \eqref{int11} and \eqref{int12} that 
\begin{align*}
\displaystyle
\frac{| \varphi - {I}_T \varphi |_{W^{m,p}({T})}}{| {\varphi} |_{W^{\ell+1,p}({T})}}
&\leq C^{A,d} C^{BH} \diam(\widehat{T})^{\ell+1-m} \{ 1+(\ell+1)  C^{\mathcal{S}} \} \left( \frac{H_T}{h_T} \right)^m h_T^{\ell+1-m}.
\end{align*}
\qed
\end{pf*}

\begin{Ex}
As the examples in \cite[Example 1.106]{ErnGue04}, we get local interpolation error estimates for a Lagrange finite element of degree $k$, a more general finite element, and the Crouzeix--Raviart finite element with $k=1$.
\begin{enumerate}
 \item For a Lagrange finite element of degree $k$, we set $V(\widehat{T}) := \mathcal{C}^0(\widehat{T})$. The condition on $\ell$ in Theorem \ref{int=thr2} is $\frac{d}{p} - 1 \< \ell \leq k$ because $W^{\ell+1,p}(\widehat{T}) \subset \mathcal{C}^0(\widehat{T})$ if $\ell+1 \> \frac{d}{p}$ according to the Sobolev imbedding theorem. 
 \item For a general finite element with $V(\widehat{T}) := \mathcal{C}^t(\widehat{T})$ and $t \in \mathbb{N}$. The condition on $\ell$ in Theorem \ref{int=thr2} is $\frac{d}{p} - 1 + t \< \ell \leq k$. When $t=1$, there is a Hermite finite element.
 \item For the Crouzeix--Raviart finite element with $k=1$, we set $V(\widehat{T}) := W^{1,1}(\widehat{T})$. The condition on $\ell$ in Theorem \ref{int=thr2} is $0 \leq \ell \leq 1$.
\end{enumerate}
\end{Ex}

\begin{Rem}
We consider optimality of the estimates. Let $T \subset \mathbb{R}^3$ be the simplex with vertices $x_1 := (0,0,0)^T$, $x_2 := (s,0,0)^T$, $x_3 := (s/2,s^{\varepsilon},0)^T$, and $x_4 := (0,0,s)^T$ ($1 \< \varepsilon \< 2$), and $0 \< s \< 1$, $s \in \mathbb{R}$. Let 
\begin{align*}
\displaystyle
\varphi(x,y,z) := x^2 + \frac{1}{4} y^2 + z^2.
\end{align*}
Let $I_T^{L}: \mathcal{C}^0(T) \to \mathcal{P}^1$ be the local Lagrange interpolation operator. We set
\begin{align*}
\displaystyle
I_T^{L} \varphi (x,y,z) := a x + by + cz + d,
\end{align*}
where $a,b,c,d \in \mathbb{R}$. For any nodes $P$ of $T$, since $I_T^{L} \varphi (P) = \varphi (P)$, we have
\begin{align*}
\displaystyle
I_T^{L} \varphi (x,y,z) = s x - \frac{1}{4} ( s^{2 - \varepsilon} - s^{\varepsilon})  y+ s z.
\end{align*}
It thus holds that
\begin{align*}
\displaystyle
(\varphi - I_T^{L} \varphi)(x,y,z)
&= x^2 + \frac{1}{4} y^2 +  z^2 - s x + \frac{1}{4} ( s^{2 - \varepsilon} - s^{\varepsilon}) y - s z,
\end{align*}
Therefore, we have
\begin{align*}
\displaystyle
\frac{|\varphi - I_T^{L} \varphi|_{W^{1,\infty}(T)}}{|\varphi|_{W^{2,\infty}(T)}}
&= \frac{\frac{1}{4}(s^{2 - \varepsilon} + s^{\varepsilon})}{2} =: I_T.
\end{align*}
By simple calculation, we have
\begin{align*}
\displaystyle
\frac{I_T}{H_T} 
= \frac{s^4 + s^{2 + 2 \varepsilon}}{48 \sqrt{2} s^3 \sqrt{(\frac{s}{2})^2 + s^{2 \varepsilon}}}
&\geq \frac{s^4 + s^{2 + 2 \varepsilon}}{24 \sqrt{10} s^4} \geq \frac{s^4}{24 \sqrt{10} s^4} = \frac{1}{24 \sqrt{10}}.
\end{align*}
We here used
\begin{align*}
\displaystyle
H_T = \frac{6 \sqrt{2} s^3 \sqrt{(\frac{s}{2})^2 + s^{2 \varepsilon}}}{s^{2 + \varepsilon}}.
\end{align*}
We conclude that
\begin{align*}
\displaystyle
|\varphi - I_T^{L} \varphi|_{W^{1,\infty}(T)}
&\geq  \frac{1}{24 \sqrt{10}} H_T |\varphi|_{W^{2,\infty}(T)}.
\end{align*}
Therefore, the parameter $H_T$ is optimal.
\end{Rem}

\section{Raviart--Thomas Interpolation Error Estimates} \label{RTint}
This section proposes error analysis for the Raviart--Thomas interpolation of arbitrary order $k \in \mathbb{N}_0$.

\subsection{Preliminaries of Error Estimates}
We first give estimates relating to the diagonal matrix \eqref{mesh1}.
\begin{lem} \label{rt=lem5}
Let $\ell$ be such that $0 \leq \ell \leq k$. It holds that, for any $\hat{v} = (\hat{v}_1,\ldots,\hat{v}_d)^T \in L^2(\widehat{T})^d$ with $\tilde{v} = (\tilde{v}_1,\ldots,\tilde{v}_d)^T := \widehat{\Psi} \hat{v}$ and $\hat{w} = (\hat{w}_1,\ldots,\hat{w}_d)^T \in H^{\ell+1}(\widehat{T})^d$ with $\tilde{w} = (\tilde{w}_1,\ldots,\tilde{w}_d)^T := \widehat{\Psi} \hat{w}$,
\begin{align}
\displaystyle
\frac{\| \tilde{v} \|_{L^2(\widetilde{T})^d}}{| \tilde{w} |_{H^{\ell+1}(\widetilde{T})^d}} 
\leq \max_{1 \leq i \leq d} \{ \alpha_i^{\ell+1} \} \frac{ \left( \sum_{i=1}^d \alpha_i^2 \| \hat{v}_i \|_{L^2(\widehat{T})}^2 \right)^{1/2}}{\left( \sum_{i=1}^d \alpha_i^2 |\hat{w}_i |_{H^{\ell+1}(\widehat{T})}^2 \right)^{1/2}}. \label{rt1}
\end{align}
\end{lem}

\begin{pf*}
From the definition of the Piola transformation, for $i=1,\ldots,d$,
\begin{align*}
\displaystyle
\tilde{w}_i(\tilde{x}) 
&= \frac{1}{|\det (\widehat{A}^{(d)})|} \sum_{j=1}^d [\widehat{A}^{(d)} ]_{ij} \hat{w}_j(\hat{x}) 
= \frac{1}{|\det (\widehat{A}^{(d)})|} \alpha_i \hat{w}_i (\hat{x}).
\end{align*}
Let  $\beta$ be a multi-index with $|\beta| = \ell + 1$. We then have
\begin{align*}
\displaystyle
\partial^{\beta} \tilde{w}_i (\tilde{x}) 
=  \frac{1}{|\det (\widehat{A}^{(d)})|} \alpha_i (\partial^{\beta} \hat{w}_i) \alpha_{1}^{- \beta_1} \cdots \alpha_{d}^{- \beta_d}.
\end{align*}
We here used $\hat{x}_j = \alpha_j^{-1} \tilde{x}_j$. 

For any $\tilde{v} \in L^{2}(\widetilde{T})^d$, from the definition of the Piola transformation, we have
\begin{align*}
\displaystyle
\| \tilde{v} \|^2_{L^2(\widetilde{T})^d}
&= \frac{1}{|\det (\widehat{A}^{(d)})|} \| \widehat{A}^{(d)} \hat{v} \|^2_{L^2(\widetilde{T})^d} 
= \frac{1}{|\det (\widehat{A}^{(d)})|} \sum_{i=1}^d \alpha_i^2 \| \hat{v}_i \|_{L^2(\widehat{T})}^2.
\end{align*}
Meanwhile, we have, for any $\tilde{w} \in H^{\ell+1}(\widetilde{T})^d$,
\begin{align*}
\displaystyle
|\tilde{w}|_{H^{\ell+1}(\widetilde{T})^d}^2 
&= \sum_{i=1}^d |\tilde{w}_i|_{H^{\ell+1}(\widetilde{T})}^2 
=  \sum_{i=1}^d \sum_{|\beta| = \ell+1}  \| \partial^{\beta} \tilde{w}_i \|^2_{L^2(\widetilde{T})} \\
&=  \frac{1}{|\det (\widehat{A}^{(d)})|} \sum_{i=1}^d \alpha_i^2 \sum_{|\beta| = \ell+1} (\alpha_{1}^{- \beta_1} \cdots \alpha_{d}^{- \beta_d})^2 \| \partial^{\beta} \hat{w}_i \|^2_{L^2(\widehat{T})} \\
&\geq \frac{1}{|\det (\widehat{A}^{(d)})|} \min_{1 \leq j \leq d} \{ \alpha_j^{-2|\beta|} \} \sum_{i=1}^d \alpha_i^2 \sum_{|\beta| = \ell+1}  \| \partial^{\beta} \hat{w}_i\|^2_{L^2(\widehat{T})}.
\end{align*}
These inequalities conclude \eqref{rt1}. 
\qed
\end{pf*}

We next give estimates relating to the matrices \eqref{mesh2} and \eqref{mesh3}.
\begin{lem}  \label{rt=lem6}
Let $\ell$ be such that $0 \leq \ell \leq k$. For any $\tilde{v} \in L^2(\widetilde{T})^d$ with ${v} := \widetilde{\Psi} \tilde{v}$ and $\tilde{w} \in H^{\ell+1}(\widetilde{T})^d$ with ${w} := \widetilde{\Psi} \tilde{w}$, we have
\begin{align}
\displaystyle
\frac{\| {v} \|_{L^2({T})^d}}{| {w} |_{H^{\ell+1}({T})^d}} &\leq C^{P,d} \frac{H_T}{h_T} \frac{\| \tilde{v} \|_{L^2(\widetilde{T})^d}}{| \tilde{w} |_{H^{\ell+1}(\widetilde{T})^d}}, \label{rt2}
\end{align}
where $C^{P,2} := 2^{\frac{\ell+1}{2}} C^{vec}$, and $C^{P,3} := \frac{2^{\ell+2}}{3} C^{vec}$, where $C^{vec}$ is a constant independent of $T$ and $\widetilde{T}$.
\end{lem}

\begin{pf*}
Using the standard estimates in \cite[Lemma 1.113]{ErnGue04}, we easily get
\begin{align}
\displaystyle
\frac{\| {v} \|_{L^2({T})^d}}{| {w} |_{H^{\ell+1}({T})^d}} &\leq C^{vec} \left( \| \widetilde{A} \|_2 \| \widetilde{A}^{-1} \|_2 \right) \| \widetilde{A} \|_2^{\ell+1}  \frac{\| \tilde{v} \|_{L^2(\widetilde{T})^d}}{| \tilde{w} |_{H^{\ell+1}(\widetilde{T})^d}}, \quad d=2,3. \label{rt3}
\end{align}
Therefore, \eqref{rt2} follows from \eqref{rt3}, \eqref{int4}, and \eqref{int5} if $d=2$ and from \eqref{rt3}, \eqref{int6}, and \eqref{int7} if $d=3$.
\qed
\end{pf*}

\subsection{Component-wise Stability of the Raviart--Thomas Interpolation on the Reference Element}
This subsection introduces the component-wise stability for the Raviart--Thomas interpolation of any order of functions in $H^1(\widehat{T})^d$. To this end, we follow \cite{AcoApe10}; see also \cite{AcoDur99}. 

We first introduce component-wise stability estimates in the reference element $\widehat{T} = \conv \{ 0,e_1, \ldots,e_d \}$. Here, $e_1, \ldots, e_d \in \mathbb{R}^d$ are the canonical basis.

\begin{lem} \label{rt=lem7}
For $k \in \mathbb{N}_0$, there exists a constant $C_1^{(i)}(k)$, $i=1,\dots,d$ such that, for all $\hat{u} = (\hat{u}_1, \ldots, \hat{u}_d)^T \in H^1(\widehat{T})^d$, 
\begin{align}
\displaystyle
\| (I_{\widehat{T}}^{RT} \hat{u})_i \|_{L^2(\widehat{T})} \leq C_1^{(i)}(k) \left( \| \hat{u}_i \|_{H^1(\widehat{T})} + \| \div \hat{u} \|_{L^2(\widehat{T})} \right), \quad i=1, \ldots, d.\label{rt4}
\end{align}
\end{lem}

\begin{pf*}
The proof is given in \cite[Lemma 3.3]{AcoApe10} for the case $d=3$. The estimate in the case $d=2$ can be proved analogously.
\qed
\end{pf*}

We next give component-wise stability estimates in the reference element $\widehat{T} = \conv \{ 0,e_1, e_1 + e_2 , e_3 \}$. 

\begin{lem} \label{rt=lem9}
For $k \in \mathbb{N}_0$, there exists a constant $C_2^{(i)}(k)$, $i=1,2,3$ such that, for all $\hat{u} = (\hat{u}_1,\hat{u}_2, \hat{u}_3)^T \in H^1(\widehat{T})^3$, 
\begin{align}
\displaystyle
\| (I_{\widehat{T}}^{RT} \hat{u})_i \|_{L^2(\widehat{T})} 
&\leq C_2^{(i)}(k) \left( \| \hat{u}_i \|_{H^1(\widehat{T})} + \sum_{j=1, j\neq i}^3 \left\| \frac{\partial \hat{u}_j}{\partial \hat{x}_j} \right\|_{L^2(\widehat{T})} \right), \quad i=1,2,3. \label{rt5}
\end{align}
\end{lem}

\begin{pf*}
The proof is given in \cite[Lemma 4.3]{AcoApe10}. We remark that our reference element in this case is different from that in  \cite[Lemma 4.3]{AcoApe10}. However, the proof can be made by analogous argument.
\qed
\end{pf*}

\subsection{Raviart--Thomas Interpolation Error Estimates}

\begin{thr}  \label{rt=thr3}
For $k \in \mathbb{N}_0$, let $\{ {T} , RT^k({T}) , {\Sigma} \}$ be the Raviart--Thomas finite element and $I_T^{RT}$ the local interpolation operator defined in \eqref{mesh9}. Let $\ell$ be such that $0 \leq \ell \leq k$. We then have the estimates
\begin{align}
\displaystyle
\| I_T^{RT} v - v \|_{L^2(T)^d} &\leq  C_I^{RT} H_T h_T^{\ell} |v|_{H^{\ell+1}(T)^d} \quad \forall v \in H^{\ell+1}(T)^d. \label{rt13}
\end{align}
Here, $C_I^{RT} := C^{P,d} C(\widehat{T},d,\ell,k) $ is a positive constant independent of $T$ while $C^{P,d}$ is the constant appearing in Lemma \ref{rt=lem6}.
\end{thr}

\begin{pf*}
Let $\hat{v} \in H^{\ell+1}(\widehat{T})$. Let ${I}_{\widetilde{T}}^{RT}$ and ${I}_{\widehat{T}}^{RT}$ respectively be the local interpolation operators on $\widetilde{T}$ and $\widehat{T}$ defined by \eqref{mesh8}, \eqref{pre8}, and \eqref{pre9}. From \eqref{rt2}, we have
\begin{align}
\displaystyle
\frac{\| I_T^{RT} v - v \|_{L^2({T})^d}}{| {v} |_{H^{\ell+1}({T})^d}} &\leq C^{P,d} \frac{H_T}{h_T} \frac{\| I_{\widetilde{T}}^{RT} \tilde{v} - \tilde{v} \|_{L^2(\widetilde{T})^d}}{| \tilde{v} |_{H^{\ell+1}(\widetilde{T})^d}}. \label{rt14}
\end{align}
If $\tilde{q} \in \mathcal{P}^{\ell} (\widetilde{T})^d \subset RT^k(\widetilde{T})$, we have $I_{\widetilde{T}}^{RT} \tilde{q} = \tilde{q}$. We therefore obtain, for any $\tilde{q} \in \mathcal{P}^{\ell} (\widetilde{T})^d$, 
\begin{align*}
\displaystyle
 \frac{\| I_{\widetilde{T}}^{RT} \tilde{v} - \tilde{v} \|_{L^2(\widetilde{T})^d}}{| \tilde{v} |_{H^{\ell+1}(\widetilde{T})^d}}
 &\leq \frac{\| I_{\widetilde{T}}^{RT} ( \tilde{v} - \tilde{q} ) \|_{L^2(\widetilde{T})^d}}{| \tilde{v} |_{H^{\ell+1}(\widetilde{T})^d}} +  \frac{\| \tilde{q} - \tilde{v} \|_{L^2(\widetilde{T})^d}}{| \tilde{v} |_{H^{\ell+1}(\widetilde{T})^d}}.
\end{align*}

Using inequality \eqref{rt1}, the component-wise stability \eqref{rt4} and \eqref{rt5}, we have
\begin{align*}
\displaystyle
&\frac{\| I_{\widetilde{T}}^{RT} ( \tilde{v} - \tilde{q} ) \|_{L^2(\widetilde{T})^d}}{| \tilde{v} |_{H^{\ell+1}(\widetilde{T})^d}} \\
&\leq h_T^{\ell+1} \frac{\left( \sum_{i=1}^d \alpha_i^2 \| \{ I_{\widehat{T}}^{RT} ( \hat{v} - \hat{q} ) \}_i \|_{L^2(\widehat{T})}^2 \right)^{1/2}}{\left( \sum_{i=1}^d \alpha_i^2 |\hat{v}_i |_{H^{\ell+1}(\widehat{T})}^2 \right)^{1/2}} \\
&\leq  C(k) h_T^{\ell+1} \frac{\left( \sum_{i=1}^d \alpha_i^2 \left\{   \| ( \hat{v} - \hat{q} )_i \|_{H^1(\widehat{T})}^2 +\sum_{j=1}^d \left\| \frac{\partial (\hat{v} - \hat{q})_j}{\partial \hat{x}_j} \right\|_{L^2(\widehat{T})}^2 \right\} \right)^{1/2}}{\left( \sum_{i=1}^d \alpha_i^2 |\hat{v}_i |_{H^{\ell+1}(\widehat{T})}^2 \right)^{1/2}}.
\end{align*}
By applying Theorem \ref{int=thr1}, we obtain
\begin{align*}
\displaystyle
&C(k) h_T^{\ell+1} \inf_{ \hat{q} \in \mathcal{P}^{\ell} (\widehat{T})^d} \frac{\left( \sum_{i=1}^d \alpha_i^2 \left\{   \| ( \hat{v} - \hat{q} )_i \|_{H^1(\widehat{T})}^2 +\sum_{j=1}^d \left\| \frac{\partial (\hat{v} - \hat{q})_j}{\partial \hat{x}_j} \right\|_{L^2(\widehat{T})}^2 \right\} \right)^{1/2}}{\left( \sum_{i=1}^d \alpha_i^2 |\hat{v}_i |_{H^{\ell+1}(\widehat{T})}^2 \right)^{1/2}} \\
&\leq C(\widehat{T},d,\ell,k)  h_T^{\ell+1} \frac{\left( \sum_{i=1}^d \alpha_i^2  |\hat{v}_i |_{H^{\ell+1}(\widehat{T})}^2 \right)^{1/2}}{\left( \sum_{i=1}^d \alpha_i^2 |\hat{v}_i |_{H^{\ell+1}(\widehat{T})}^2 \right)^{1/2}}.
\end{align*}

Furthermore, using inequality \eqref{rt1}, we have
\begin{align*}
\displaystyle
  \frac{\| \tilde{q} - \tilde{v} \|_{L^2(\widetilde{T})^d}}{| \tilde{v} |_{H^{\ell+1}(\widetilde{T})^d}}
  &\leq h_T^{\ell+1} \frac{ \left( \sum_{i=1}^d \alpha_i^2 \| \hat{q}_i - \hat{v}_i \|_{L^2(\widehat{T})}^2 \right)^{1/2}}{\left( \sum_{i=1}^d \alpha_i^2 |\hat{v}_i |_{H^{\ell+1}(\widehat{T})}^2 \right)^{1/2}}.
\end{align*}
We apply Theorem \ref{int=thr1} to obtain
\begin{align*}
\displaystyle
&h_T^{\ell+1} \inf_{ \hat{q} \in \mathcal{P}^{\ell} (\widehat{T})^d} \frac{ \left( \sum_{i=1}^d \alpha_i^2 \| \hat{q}_i - \hat{v}_i \|_{L^2(\widehat{T})}^2 \right)^{1/2}}{\left( \sum_{i=1}^d \alpha_i^2 |\hat{v}_i |_{H^{\ell+1}(\widehat{T})}^2 \right)^{1/2}} \\
&\quad \leq C(\widehat{T},d,\ell)  h_T^{\ell+1} \frac{\left( \sum_{i=1}^d \alpha_i^2  |\hat{v}_i |_{H^{\ell+1}(\widehat{T})}^2 \right)^{1/2}}{\left( \sum_{i=1}^d \alpha_i^2 |\hat{v}_i |_{H^{\ell+1}(\widehat{T})}^2 \right)^{1/2}}.
\end{align*}

We thus have
\begin{align}
\displaystyle
 \frac{\| I_{\widetilde{T}}^{RT} \tilde{v} - \tilde{v} \|_{L^2(\widetilde{T})^d}}{| \tilde{v} |_{H^{\ell+1}(\widetilde{T})^d}}
 \leq  C(\widehat{T},d,\ell,k)  h_T^{\ell+1}. \label{rt15}
\end{align}
We conclude from \eqref{rt14} and \eqref{rt15} that
\begin{align*}
\displaystyle
\frac{\| I_T^{RT} v - v \|_{L^2({T})^d}}{| {v} |_{H^{\ell+1}({T})^d}} &\leq C^{P,d}C(\widehat{T},d,\ell.k)  h_T^{\ell+1}.
\end{align*}
\qed
\end{pf*}

%
%

\begin{acknowledgements}
This work was supported by JSPS KAKENHI Grant Number \\ JP16H03950. We would like to thank the anonymous referee for the valuable comments. 
\end{acknowledgements}

%
%



\end{document}